\documentclass[12pt, reqno, a4paper]{amsart}

\usepackage{amsthm, amsmath, amssymb, amscd,latexsym}
\usepackage[pdftex]{graphicx, color}
\usepackage{enumerate}

\def\C{\mathbb C}

\def\D{\mathbb D}

\def\N{\mathbb N}

\def\wt{\widetilde}

\def\noin{\noindent}

\def\vep{\varepsilon}

\theoremstyle{plain}
\newtheorem{thm}{Theorem}[section]

\newtheorem{lem}[thm]{Lemma}
\newtheorem*{thmA*}{Theorem A}
\newtheorem*{thmA'*}{Theorem A'}
\newtheorem*{corA'*}{Corollary A'}
\newtheorem*{thmB*}{Theorem B}
\newtheorem*{thmC*}{Theorem C}
\newtheorem*{corC*}{Corollary C}
\newtheorem*{corD*}{Corollary D}
\newtheorem*{thmD*}{Theorem D}
\newtheorem*{thmE*}{Theorem E}
\newtheorem*{thmF*}{Theorem F}
\newtheorem*{conj*}{Conjecture}
\newtheorem*{thm*}{Theorem}
\theoremstyle{definition}

\newtheorem{Question*}[thm]{Question}
\newtheorem*{defn*}{Definition}
\newtheorem*{rem*}{Remark}

\makeatletter
\@addtoreset{equation}{section} 
\makeatother

\newcommand{\bs}[1]{\boldsymbol{#1}}
\renewcommand{\Bar}{\overline}

\newcommand{\abs}[1]{{\left| #1 \right|}}
\newcommand{\paren}[1]{{\left( #1 \right)}}
\newcommand{\braces}[1]{{\left\{ #1 \right\}}}
\newcommand{\al}{\alpha}
\newcommand{\e}{\varepsilon}
\newcommand{\lam}{\lambda}
\newcommand{\cM}{{\mathcal{M}}}

\newcommand{\cc}{\circ}
\newcommand{\st}{\ | \ }
\newcommand{\QED}{\hfill $\blacksquare$}

\newcommand{\sminus}{\smallsetminus}

\newcommand{\parag}[1]{\medskip \par\noindent{\bf #1}}

\begin{document}

\title[]{Julia sets appear quasiconformally
 in the Mandelbrot set, II: A parabolic proof}

\author[]{Tomoki Kawahira and Masashi Kisaka}
\setcounter{footnote}{-1}

\address{
Graduate School of Economics\\
Hitotsubashi University\\
Tokyo 186-8601 \\
Japan}
\email{t.kawahira@r.hit-u.ac.jp}

\address{
Department of Mathematical Sciences \\
Graduate School of Human and Environmental Studies \\
Kyoto University \\
Kyoto 606-8501 \\
Japan }
\email{kisaka@math.h.kyoto-u.ac.jp}

\begin{abstract}
Following the ideas of A.~Douady, 
we give an alternative
proof of the authors' result: 
for any boundary point $c_0$ of the Mandelbrot set $M$,
we can find small quasiconformal copies of $M$ in $M$
that are encaged in nested quasiconformal copies of 
the totally disconnected Julia set 
of a parameter arbitrarily close to $c_0$.
\end{abstract}

\maketitle

\section{Introduction}
In the paper \cite{Douady 2000} by A.~Douady, X.~Buff, R.~Devaney,
and P.~Sentenac titled ``Baby Mandelbrot sets are born
in cauliflowers," 
they showed that we can find small quasiconformal copies of 
the Mandelbrot set $M$ in $M$ that are encaged in nested quasiconformal copies of 
an imploded cauliflower 
(the Julia set of $z \mapsto z^2+1/4+\e$ for small $\e >0$). 
Indeed, we can always find such a copy near 
the cusp point of the primitive small copies 
of $M$ and we can visually observe imploded and nested 
cauliflowers around them.
The proof relies on the parabolic implosion technique
developed by Douady, Lavaurs, and Shishikura.

Later in \cite{Kawahira-Kisaka 2018}, 
the authors extended this result 
and showed that  
fairly large varieties of quadratic Julia sets 
appear in $M$, 
but the proof presented in that paper is based on 
the shooting technique around Misiurewicz parameters.
The aim of this paper is to present an alternative proof 
{\it \`a la Douady},
replacing ``Misiurewicz" by ``parabolic."

\parag{The main result.}
We will loosely follow Douady's original 
notation in \cite{Douady 2000}. 
We set 
\begin{eqnarray*}
& &  D(R) := \{ z \in \C \ | \ |z| < R \}, \quad
\D:=D(1), \quad
D(\alpha, R) := \{ z \in \C \ | \ |z-\alpha| < R \}, \\
& & A(r, R) := \{ z \in \C \ | \ r < |z| < R \} \quad (0 < r < R).
\end{eqnarray*}
For the quadratic map 
$P_c(z) := z^2 +c \ (c \in \C)$,
let $K(P_c)$ and $J(P_c)$ denote the 
filled Julia set and the Julia set respectively.
Now we choose any $\sigma \in \C \sminus M$
such that $J(P_{\sigma})$ is a Cantor set. 
We also choose an $R>1$ such that
$$
J(P_{\sigma}) 
\subset 
A(R^{-1/2} , R^{1/2}),
$$
and define the {\it rescaled Julia set}
$\Gamma_0(\sigma)$
by
$$
\Gamma_0(\sigma) 
:= 
J(P_{\sigma}) \times R^{3/2}
=
\braces{ R^{3/2} \,z \ \big| \ z \in J(P_{\sigma}) }
$$ 
in such a way that 
$\Gamma_0(\sigma)$ is contained in $A(R, R^2)$.

Let $\Gamma_m(\sigma) \ (m \in \N)$ be the inverse image of $\Gamma_0(\sigma)$ by 
$z \mapsto z^{2^m}$. 
Then the sets 
$\Gamma_0(\sigma),\, \Gamma_1(\sigma),\, \Gamma_2(\sigma), \ldots$ 
are mutually disjoint since 
$
\Gamma_m(\sigma) \subset A(R^{2^{-m}}, R^{2^{-m+1}}).
$

We define the {\it decorated Mandelbrot set} 
$\mathcal{M}(\sigma)$
by 
$$
\mathcal{M}(\sigma) 
:= M \cup \Phi_M^{-1}
\Big( \bigcup_{m=0}^\infty \Gamma_m(\sigma) \Big),
$$
where $\Phi_M : \C \smallsetminus M \to \C \smallsetminus \overline{\D}$ is 
the conformal isomorphism with 
$\Phi_M(c)/c \to 1 \ \text{as} \ |c| \to \infty$.

\fboxsep=0pt
\fboxrule=1pt
\begin{figure}[bp]
\begin{center}
(i)\\[.5em]
\fbox{\includegraphics[width=.18\textwidth, bb = 0 0 1000 1002]{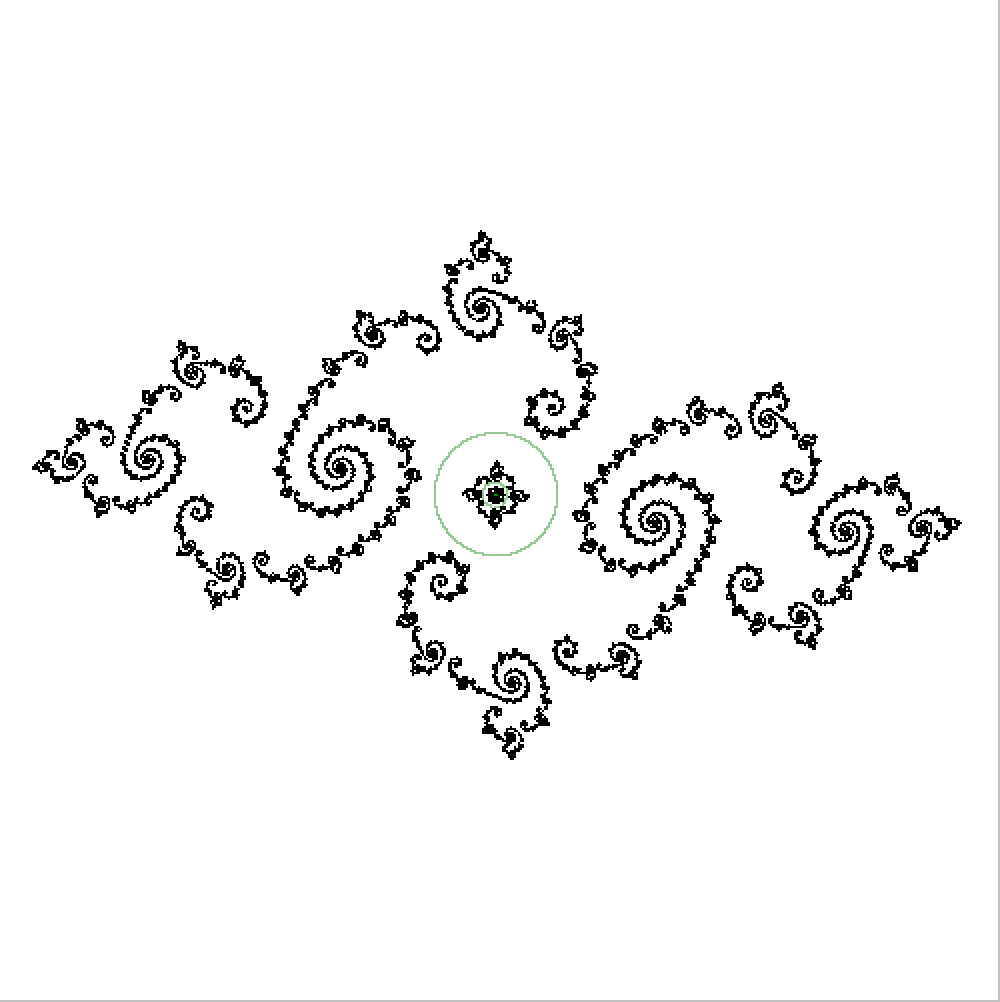}}
\fbox{\includegraphics[width=.18\textwidth, bb = 0 0 1002 1000]{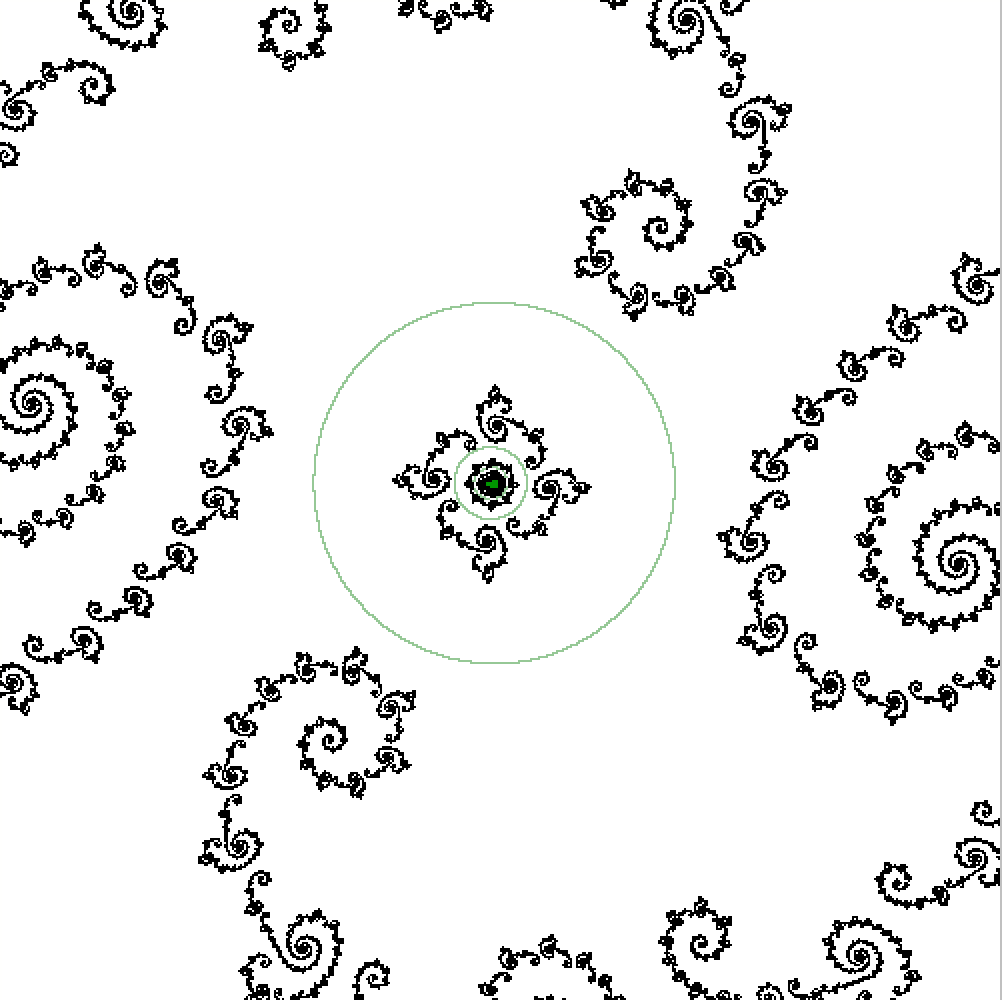}}
\fbox{\includegraphics[width=.18\textwidth, bb = 0 0 1002 998]{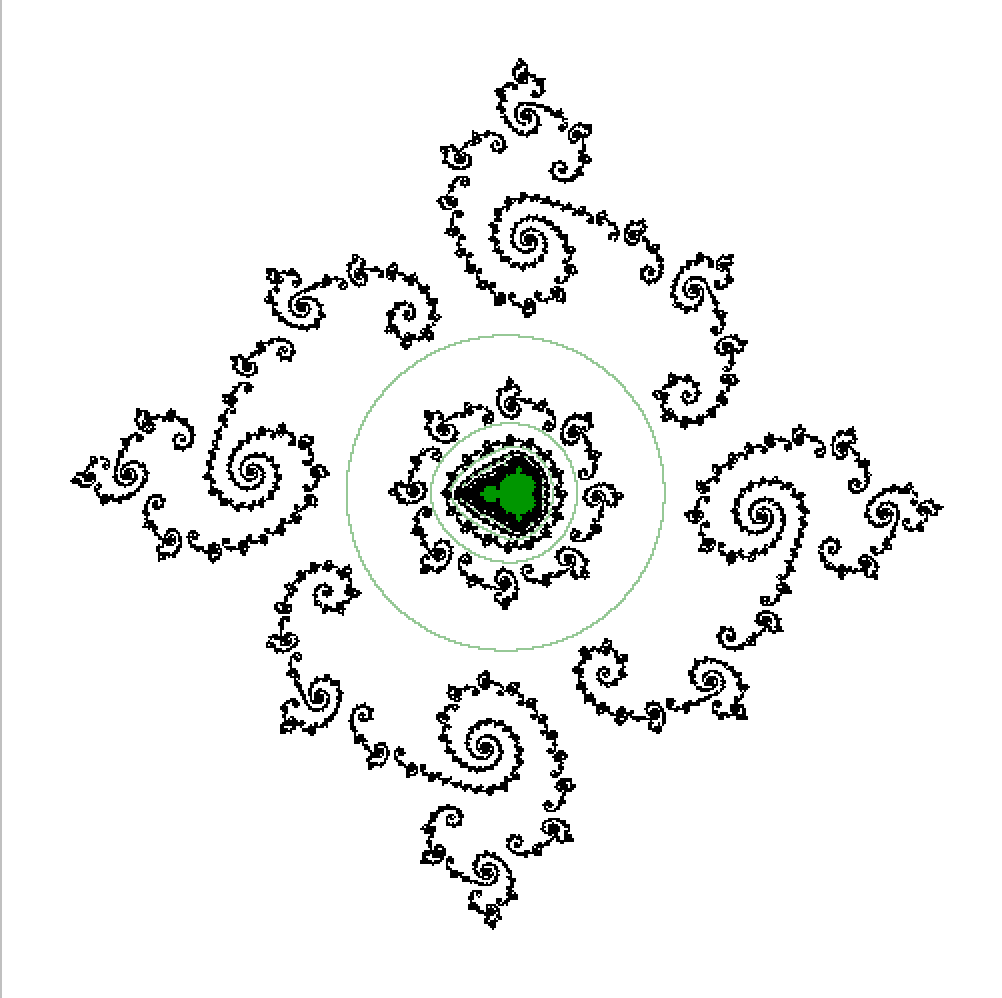}}
\fbox{\includegraphics[width=.18\textwidth, bb = 0 0 1002 1002]{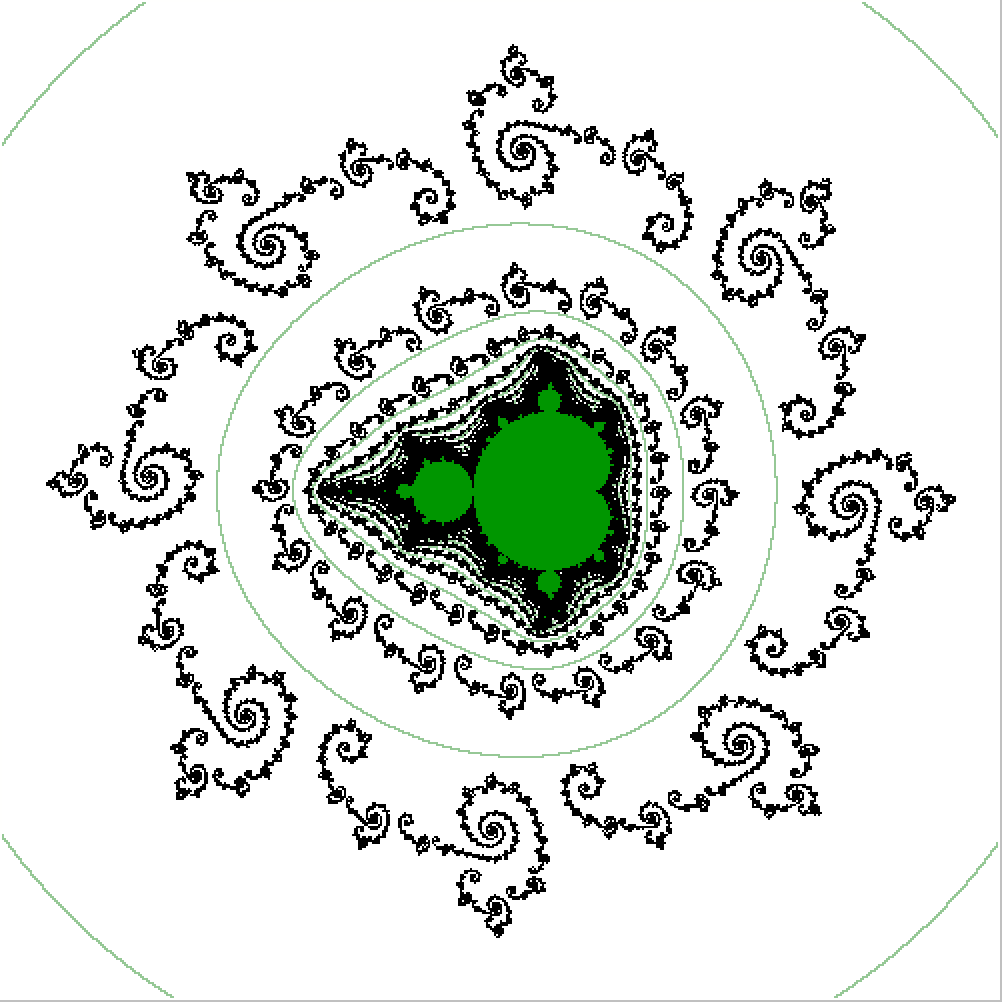}}
\fbox{\includegraphics[width=.18\textwidth, bb = 0 0 1002 1002]{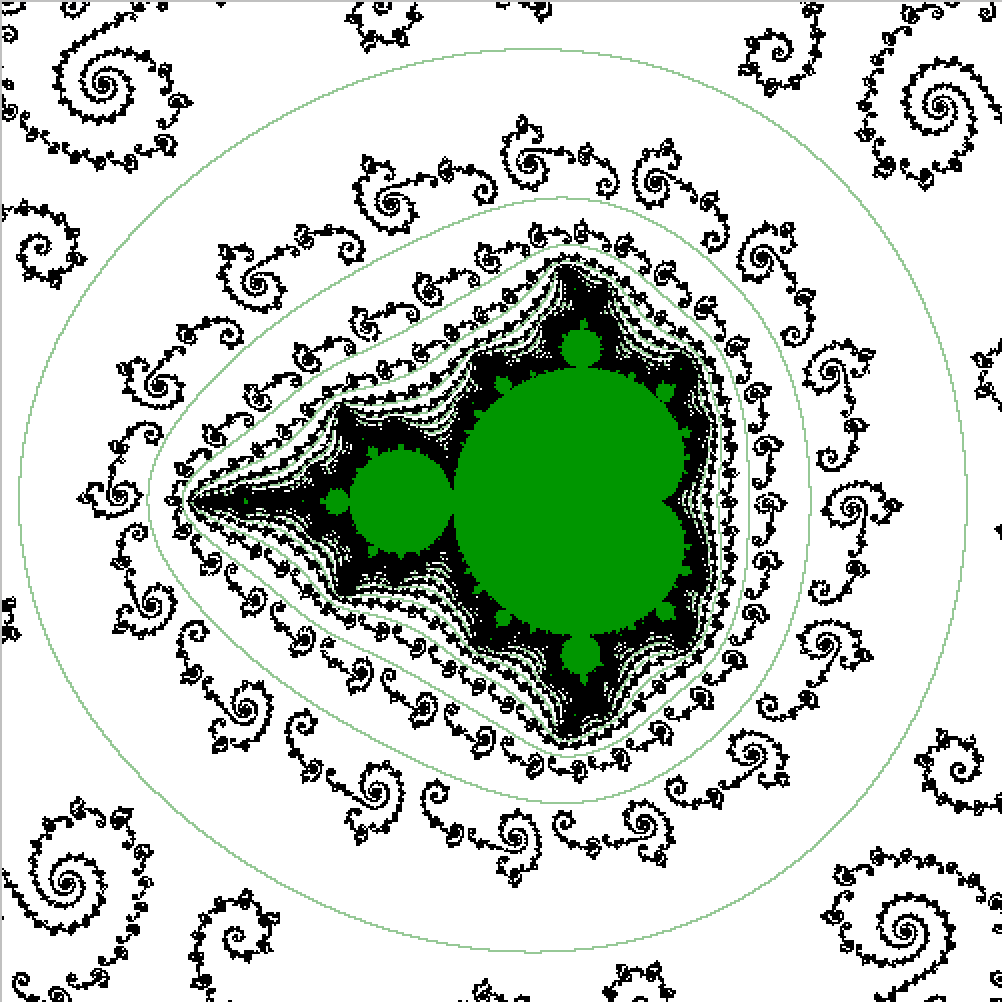}}
\\[.5em]
(ii) \\[.5em]
\fbox{\includegraphics[width=.18\textwidth, bb = 0 0 500 500]{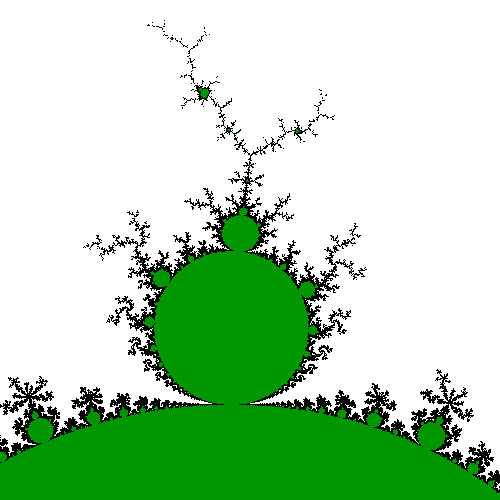}}
\fbox{\includegraphics[width=.18\textwidth, bb = 0 0 500 500]{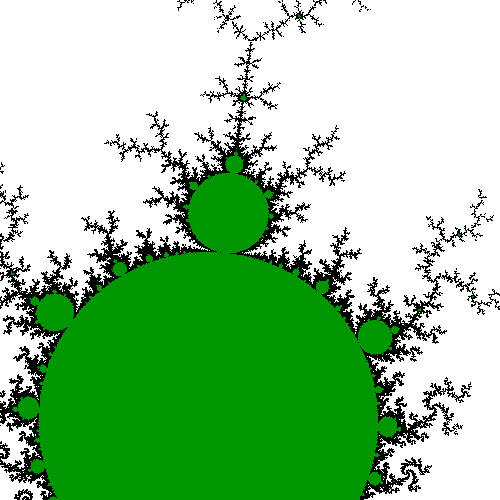}}
\fbox{\includegraphics[width=.18\textwidth, bb = 0 0 500 500]{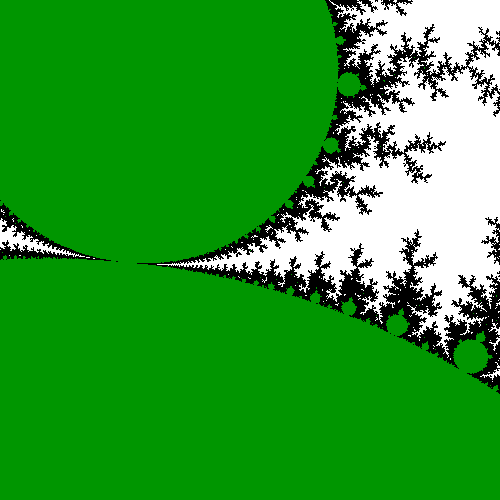}}
\fbox{\includegraphics[width=.18\textwidth, bb = 0 0 500 500]{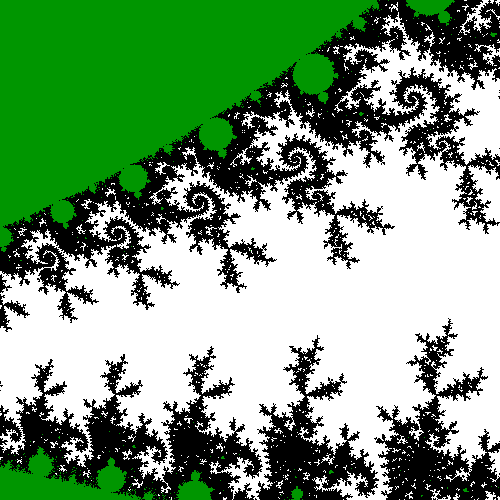}}
\fbox{\includegraphics[width=.18\textwidth, bb = 0 0 500 500]{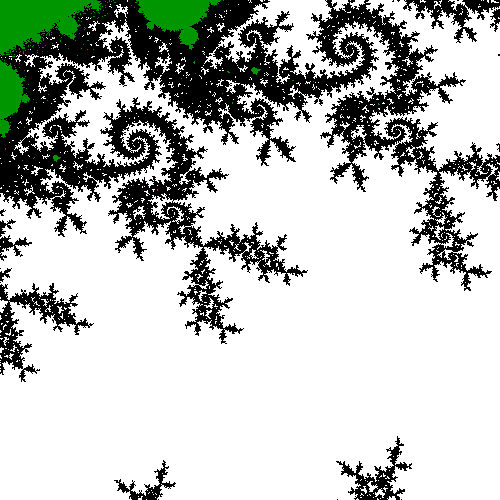}}
\\[.5em]
\fbox{\includegraphics[width=.18\textwidth, bb = 0 0 500 500]{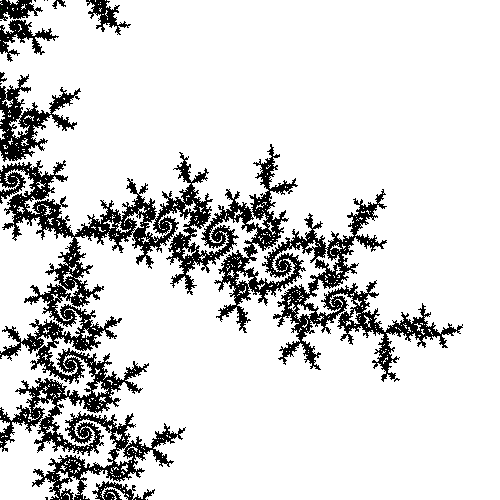}}
\fbox{\includegraphics[width=.18\textwidth, bb = 0 0 500 500]{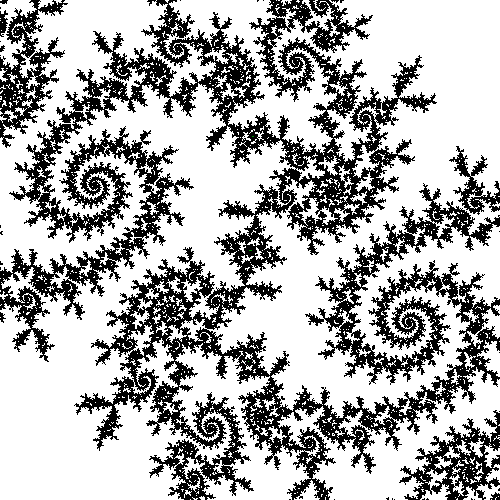}}
\fbox{\includegraphics[width=.18\textwidth, bb = 0 0 500 500]{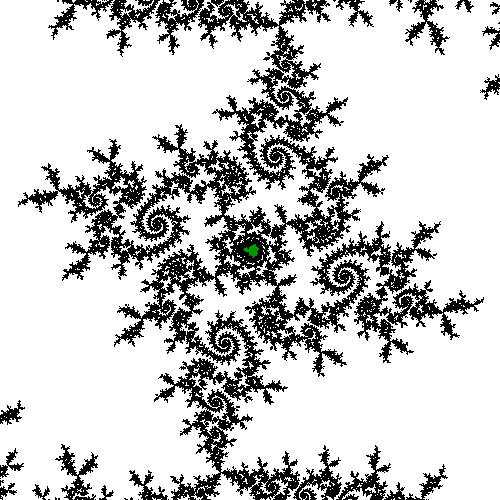}}
\fbox{\includegraphics[width=.18\textwidth, bb = 0 0 500 500]{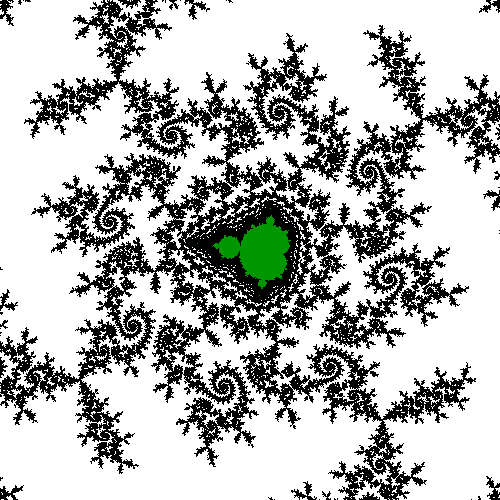}}
\fbox{\includegraphics[width=.18\textwidth, bb = 0 0 500 500]{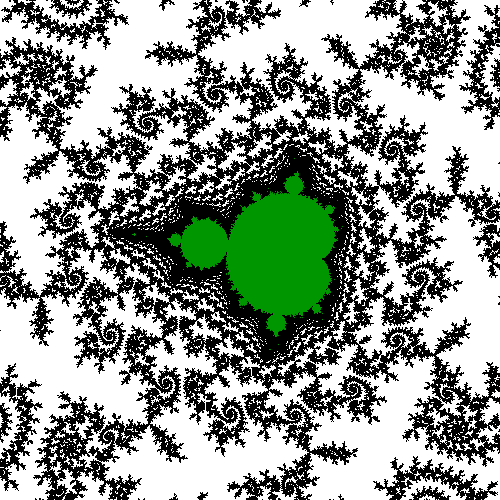}}
\\[.5em]
(iii) \\[.5em]
\fbox{\includegraphics[width=.18\textwidth, bb = 0 0 500 500]{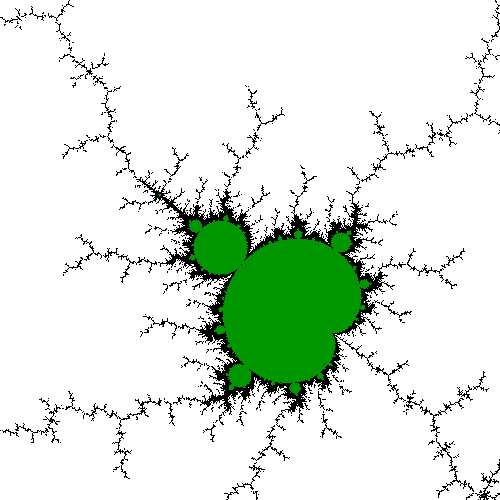}}
\fbox{\includegraphics[width=.18\textwidth, bb = 0 0 500 500]{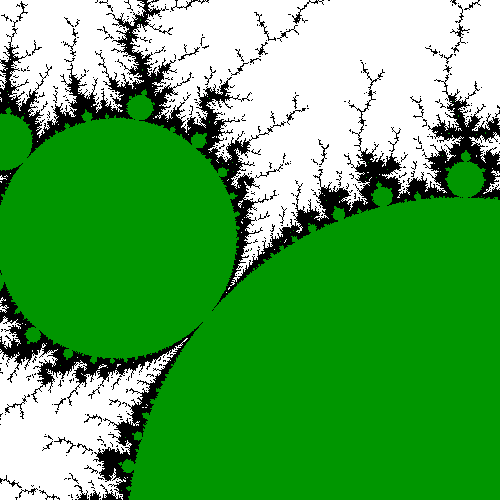}}
\fbox{\includegraphics[width=.18\textwidth, bb = 0 0 500 500]{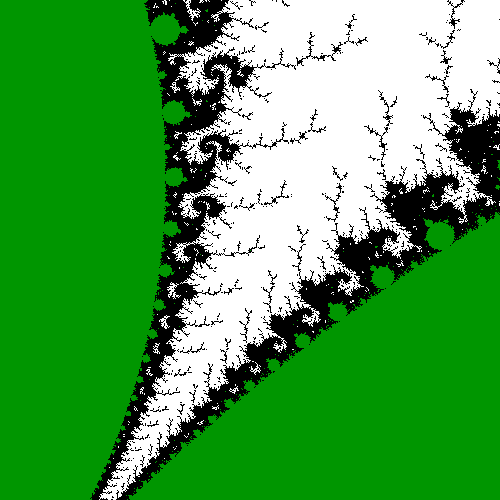}}
\fbox{\includegraphics[width=.18\textwidth, bb = 0 0 500 500]{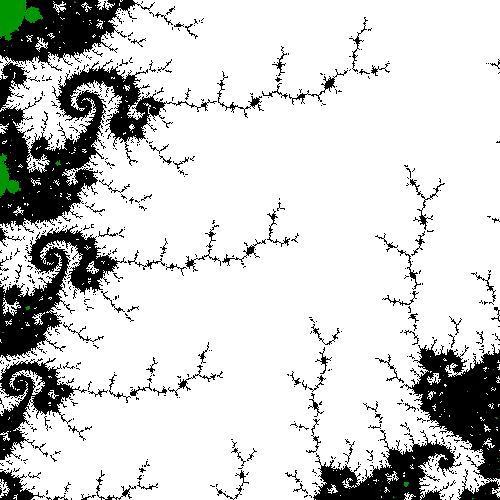}}
\fbox{\includegraphics[width=.18\textwidth, bb = 0 0 500 500]{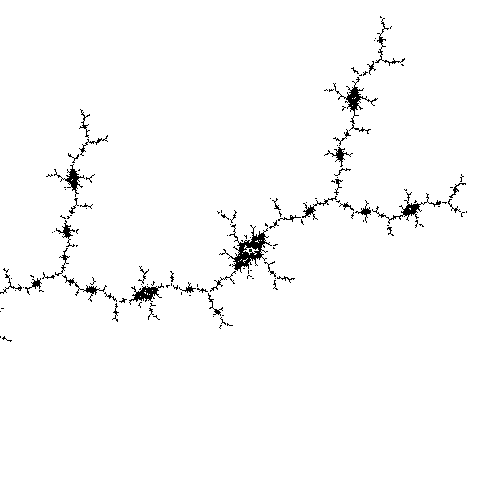}}
\\[.5em]
\fbox{\includegraphics[width=.18\textwidth, bb = 0 0 500 500]{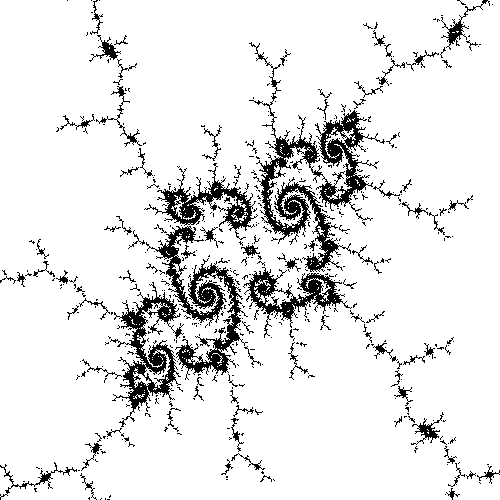}}
\fbox{\includegraphics[width=.18\textwidth, bb = 0 0 500 500]{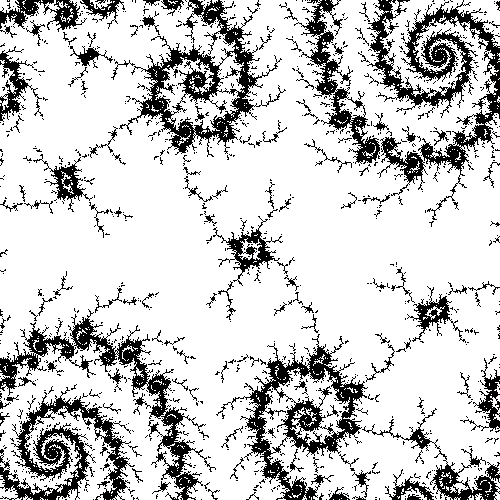}}
\fbox{\includegraphics[width=.18\textwidth, bb = 0 0 500 500]{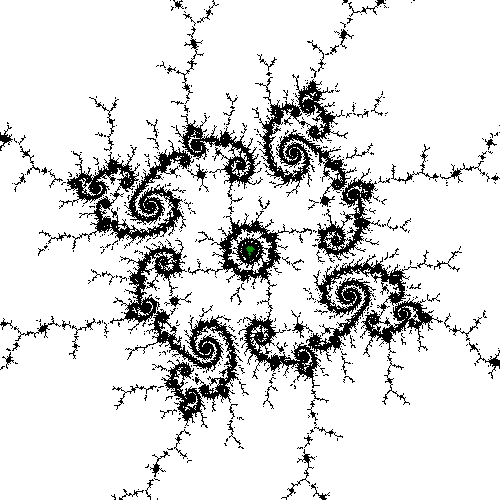}}
\fbox{\includegraphics[width=.18\textwidth, bb = 0 0 500 500]{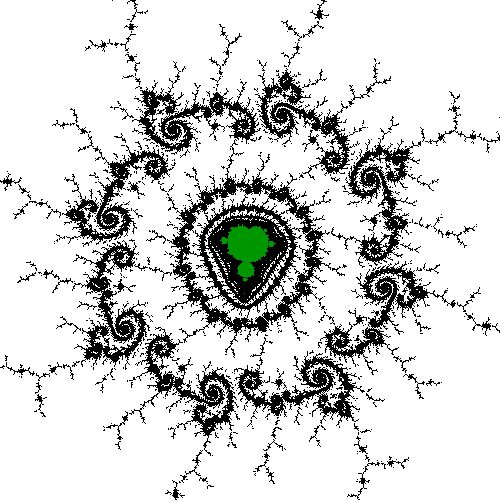}}
\fbox{\includegraphics[width=.18\textwidth, bb = 0 0 500 500]{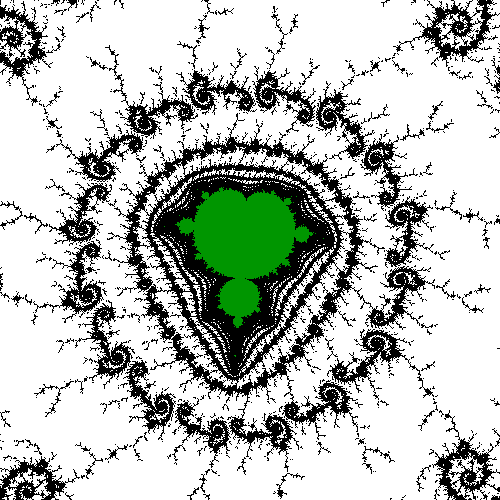}}
\end{center}
\caption{\small
(i): The decorated Mandelbrot set $\mathcal{M}(\sigma)$ for $\sigma=-0.77+0.18 i$
(close to the parabolic parameter $c_0=-0.75$). 
(ii) and (iii): Embedded quasiconformal copies of 
$\mathcal{M}(\sigma)$ 
above near satellite and primitive small Mandelbrot sets.
}
\label{fig_satellite_primitive}
\end{figure}

Let $X$ and $Y$ be non-empty compact sets in $\C$. 
We say 
{\it $Y$ contains a 
quasiconformal copy of $X$} 
if there is a quasiconformal map 
$\chi$ on a neighborhood of $X$ such that
$\chi(X) \subset Y$ and $\chi(\partial X) \subset \partial Y$. 
Following A.~Douady, we also say 
{\it $X$ appears quasiconformally in $Y$.}
Note that
the condition $\chi(\partial X) \subset \partial Y$ is to exclude the case
$\chi(X) \subset \text{int}(Y)$.

Now we are ready to state the main theorem of this paper:

\begin{thm}[\bf Julia sets appear quasiconformally]
\label{thm_main}
For any choices of 
$c_0 \in \partial M$ and $\e>0$, 
there exists a parameter 
$$
\sigma \in \D(c_0,\e) \sminus M
$$
such that $M$ contains a quasiconformal
copy of the decorated Mandelbrot set $\cM(\sigma)$. 
Moreover, one can find such a copy in any open disk intersecting with $\partial M$.
\end{thm}
\noindent
Since $\cM(\sigma)$ contains a rescaled Julia set 
$\Gamma_0(\sigma)=J(P_{\sigma}) \times R^{3/2}$,
we may say that 
{\it the Julia set $J(P_{\sigma})$ appears quasiconformally in $M$.} 

Note that if $K(P_{c_0})$ has empty interior
(i.e., $P_{c_0}$ has no parabolic basins nor Siegel disks), 
then $J(P_{\sigma})$ tends to $J(P_{c_0})$
in the Hausdorff topology as $\sigma \to c_0$.
Even in the case when the interior of $K(P_{c_0})$ is non-empty,
$J(P_{\sigma})$ is contained in the $\eta$-neighborhood of 
$K(P_{c_0})$, and the $\eta$-neighborhood of 
$J(P_{\sigma})$ contains $J(P_{c_0})$
for any given $\eta >0$ if $\sigma$ is sufficiently close to $c_0$. See \cite{Douady 1994}.
This explains why we can find structures
that resemble the Julia set $J(P_{c_0})$ everywhere 
in the boundary of the Mandelbrot set.

\parag{Small Mandelbrot sets.}
The statement of the authors' original theorem 
(Theorem A of \cite{Kawahira-Kisaka 2018}) 
has more information about the location of the copy.
The precise version can be described in terms 
of {\it tuning} and {\it small Mandelbrot sets}.

Let $s_0 \ne 0$ be a superattracting parameter
such that the period of the critical point $0$ 
is more than one.
By the Douady-Hubbard tuning theorem
\cite[Th\'eor\`eme 1 du Modulation]{Haissinsky 2000},
there exists a unique compact subset
$M_{s_0}$ of $M$ associated with a canonical homeomorphism 
$\chi_{s_0} : M_{s_0} \to M$ such that $\chi_{s_0}(s_0) = 0$. 
We also denote $M_{s_0}$ by $s_0 \perp M$
and call it the {\it small Mandelbrot set 
with center $s_0$}. 
Similarly, for $c_0 \in M$, 
let $s_0 \perp c_0$ denote the parameter
$\chi_{s_0}^{-1}(c_0)$ in $M_{s_0}$.

Theorem \ref{thm_main} can be derived from the following result:

\begin{thm}[\bf Theorem A of \cite{Kawahira-Kisaka 2018}]
\label{thm_A'}
Let $c_0$ be any parameter in $\partial M$,
and $M_{s_0}$ be any small Mandelbrot set
with center $s_0 \ne 0$.
Let $c_1:=s_0 \perp c_0 \in \partial M_{s_0}$.
Then for any $\vep > 0$ and $\vep' > 0$, 
there exists an $\eta \in \C$ with $|\eta| < \vep$ and
$c_0+\eta \notin M$ such that
${\mathcal M}(c_0+\eta)$ appears quasiconformally in 
$M \cap \Bar{D(c_1 ,\vep')}$.
In particular, the Cantor
Julia set $J(P_{c_0+\eta})$ appears quasiconformally in $M$.
\end{thm}

To obtain Theorem \ref{thm_main} from Theorem \ref{thm_A'},
 we take any $c_0 \in \partial M$
and any open disk $D$ intersecting with $\partial M$.
 Since $D$ contains a Misiurewicz parameter $c$,
and there is a sequence of small Mandelbrot sets 
converging to the $c$
(see \cite[Chapter V]{Douady-Hubbard 1985}),
we may take a small Mandelbrot set $M_{s_0}$ in $D$.
By Theorem \ref{thm_A'} we can find 
a quasiconformal copy of $\cM(\sigma)$ in $D$ 
for some $\sigma=c_0+\eta \notin M$ with arbitrarily small $|\eta|$.

\parag{Organization of the paper.}
In Section 2, we recall some fundamental facts about 
quadratic-like maps and tuning of the Mandelbrot set.
In Section 3 we summarize the parabolic implosion technique 
developed by Douady, Lavaurs, and Shishikura. 
Sections 4 to 7 are devoted for the main steps (P1)--(P4) 
of the proof of Theorem \ref{thm_A'}.
 
\parag{Notes.}
It is well-known that the Mandelbrot set 
inherits the structure of the quadratic Julia sets.
See 
\cite{Berteloot 2023}
\cite{Buff-Henriksen 2001},
\cite{Graczyk Swiatek 2011},
\cite{Graczyk Swiatek 2023},
\cite{Jung 2015},
\cite{Kawahira 2014},
\cite{MTU 1995},
\cite{MNTU 2000},
\cite{Peitgen-Saupe 1988},
\cite{Rivera-Letelier 2001},
\cite{Shishikura 1998}, and 
\cite{Tan Lei 1990} for example.
For more details, readers may consult 
the introduction of \cite{Kawahira-Kisaka 2018}.

Similar phenomena to the quadratic family 
are observed also in the 
unicritical family $\{ z^d+c \}_{c \,\in\, \C}$.
In particular, similar theorems to our main results 
should be 
formalized and proved in the same manner 
as for the quadratic case.

\section{Quadratic-like maps and renormalization}

In this section we briefly recall a fundamental theory 
of quadratic-like maps.
See \cite{Douady-Hubbard 1985}, 
\cite[Th\'eor\`eme 1 du Modulation]{Haissinsky 2000}, 
\cite{Milnor 2000}, and \cite{Lyubich Book} for more details.

\parag{Quadratic-like mappings.}
Let $U'$ and $U$ be topological disks in $\C$ 
satisfying $U' \Subset U$ (i.e., $\Bar{U'} \subset U$).
A holomorphic map 
$h:U' \to U$ is called a {\it quadratic-like map}
if $h$ is a proper branched covering of degree two. 
We define the {\it filled Julia set} $K(h)$ 
and the {\it Julia set} $J(h)$ of $h$ by
$$
K(h) := \bigcap_{n=0}^\infty h^{-n} (U'),
\quad\text{and}\quad
J(h):=\partial K(h).
$$
By the Douady-Hubbard straightening theorem
\cite[p296, Theorem 1]{Douady-Hubbard 1985},
there exists a quadratic map $P_c(z)=z^2+c$ and 
a quasiconformal map $\phi:U \to \phi(U)$
such that $\phi \cc h= P_c \cc \phi$
and $\overline{\partial} \phi =0$ a.e.~on $K(h)$.
Such a parameter $c$ is unique when $K(h)$ is connected,
and we say {\it the quadratic-like map $h$ is hybrid 
equivalent to $P_c$.}

\parag{Primitive vs.~satellite.}
Let $s_0 \ne 0$ be any superattracting parameter
(given in the statement of Theorem \ref{thm_A'}) 
such that the period of the critical point $0$ 
is exactly $p \ge 2$.
We say the small Mandelbrot set $M_{s_0}$ is {\it primitive} 
if $P_{s_0 \perp (1/4)}$ has a parabolic
periodic point with a single petal. 
Otherwise we say $M_{s_0}$ is {\it satellite},
in which case $P_{s_0 \perp (1/4)}$ has a parabolic
periodic point with more than one petal.
One can visually distinguish them
by looking at the hyperbolic component $X_0$ containing $s_0$:
It is primitive if the boundary of $X_0$
has a cuspidal point at $s_0 \perp (1/4)$;
or it is satellite if there is another hyperbolic component 
$X_1$ such that 
$\partial X_1 \cap \partial X_0=\{s_0 \perp (1/4)\}$.
(See Figure 4 of \cite{Kawahira-Kisaka 2018}, for example.)

By the Douady-Hubbard tuning theorem
\cite[p.42, Th\'eor\`eme 1 du Modulation]{Haissinsky 2000},
there exists a simply connected domain $\Lambda=\Lambda_{s_0}$ in 
the parameter plane  with the following properties:
\begin{itemize}
\item
For any $c \in \Lambda$, $P_c$ is renormalizable with period $p$. 
More precisely, there exist two Jordan domains $\widetilde{U}_c'$ and $\widetilde{U}_c$ with piecewise analytic boundaries such that
$$
  f_c := P_c^p|_{\widetilde{U}_c'} : \widetilde{U}_c' \to \widetilde{U}_c
$$
is a quadratic-like map with a critical point $0 \in \widetilde{U}_c'$. 
In particular, the boundaries of $\widetilde{U}_c'$ and 
$\widetilde{U}_c$ move holomorphically with respect to $c$ over $\Lambda$.
\item
There exists a canonical homeomorphism
$\chi_{s_0}:\Lambda \to \chi_{s_0}(\Lambda)$
such that 
$M\sminus \{1/4\} \subset \chi_{s_0}(\Lambda)$
and
for each $c \in \Lambda$, 
$f_c:\widetilde{U}_c' \to \widetilde{U}_c$ 
is hybrid equivalent to $z \mapsto z^2 +\chi_{s_0}(c)$.
\item
In both cases, 
$\chi_{s_0}^{-1}$ restricted to $M \sminus \{1/4\}$ 
extends to a homeomorphism 
$\chi_{s_0}^{-1}:M \to \chi_{s_0}^{-1}(M)$.
This image $\chi_{s_0}^{-1}(M)$ is the small Mandelbrot set $M_{s_0}$.
\item
If $M_{s_0}$ is a primitive small Mandelbrot set, then 
$M_{s_0}\subset \Lambda$.
\item
If $M_{s_0}$ is a satellite small Mandelbrot set,   
then $M_{s_0}\sminus\{s_0 \perp (1/4)\}\subset \Lambda$.
\end{itemize}
This family $\{f_c\}_{c \,\in \,\Lambda}$
of quadratic-like maps above
is a key ingredient of the proof of Theorem \ref{thm_A'}.

\parag{Outline of the proof of Theorem \ref{thm_A'}.}
First, since the parabolic parameters 
(that is, parameter $c$ for which $P_c$ has 
a parabolic periodic point) are dense in $\partial M$,
we may assume that $c_0$ in the statement is 
a parabolic parameter.
Then for the parameter $c_1:=s_0 \perp c_0$,
we consider perturbation $P_c$ of $P_{c_1}$
such that the parameter $c$ ranges over a sector $S$ 
attached to $c_1$.

Next the proof breaks into four steps (P1)--(P4):
In Steps (P1) and (P2), 
we construct two families of quadratic-like maps
$\bs{f}= \{f_c:U_c' \to U_c\}_{c \, \in \, S \cap \Lambda}$
and 
$\bs{G}= \{G_c:V_c' \to U_c\}_{c \, \in \, W}$,
that are ``nested" in both dynamical and parameter spaces
in the sense that $W \Subset S \cap \Lambda$ 
and $V_c' \Subset U_c'$ for 
$c \in W$.
Then in Steps (P3) and (P4), we check that 
the first family 
$\{f_c\}_{c \, \in \, S\cap \Lambda}$
restricted on $c \in W$ provides a stable quasiconformal copy 
of the Julia set $J(P_{c_0+\eta})$ in the statement,
and the second family
$\{G_c\}_{c \, \in \, W}$
provides a quasiconformal copy of 
the decorated Mandelbrot set $\cM(c_0+\eta)$.

\section{Fatou coordinates}
In this section we recall some fundamental facts 
about Fatou coordinates and their perturbations
near parabolic parameters that will be mainly used in 
Step (P2) (Section 5).

Agan let $M_{s_0}$ be the small Mandelbrot set with 
center $s_0 \neq 0$
such that $0$ is a periodic point of period $p \ge 2$,
and let $\Lambda=\Lambda_{s_0}$ be the simply connected domain 
where the family 
$
\{
f_c := P_c^p|_{\widetilde{U}_c'} : \widetilde{U}_c' \to \widetilde{U}_c
\}_{c \, \in \, \Lambda}
$
of quadratic-like maps (given in the previous section) 
is defined.
As we have remarked, 
since the parabolic parameters are dense in $\partial M$,
we may assume that 
$c_0 \in \partial M$ 
in Theorem \ref{thm_A'} is a parabolic parameter. 
Let $c_1:=s_0 \perp c_0 \in M_{s_0}$.
Note that $c_1$ is also a parabolic parameter.

\parag{A pair of petals and the Fatou coordinates.}
We start with the global dynamics of $P_{c_1}:\C \to \C$
including that of $f_{c_1}= P_{c_1}^p|_{\widetilde{U}_{c_1}'}$.
Let $\Delta$ be the Fatou component of $K(f_{c_1})$
(i.e., the connected component of the interior of $K(f_{c_1})$) 
containing $0$. 
The boundary $\partial \Delta$ contains
a unique parabolic periodic point $q_{c_1}$ of $f_{c_1}$ 
(resp. $P_{c_1}$)
of period $k$ 
(resp. $kp$).
The multiplier $\mu_{c_1} := (f_{c_1}^k)'(q_{c_1})$ is of the form
$$
\mu_{c_1} := 
(f_{c_1}^k)'(q_{c_1})
=e^{2 \pi i \nu'/\nu},
$$
where $\nu'$ and $\nu$ are coprime positive integers.
Since $P_{c_1}$ has only one critical point,
$q_{c_1}$ has $\nu$-petals.
That is, by choosing an appropriate local coordinate $w=\psi_{c_1}(z)$ 
near $q_{c_1}$ with $\psi_{c_1}(q_{c_1})=0$,
we have 
$$
\psi_{c_1} \circ f_{c_1}^{k\nu} \circ \psi_{c_1}^{-1}(w)=w(1 + w^{\nu}+O(w^{2\nu})).
$$
See \cite[Proof of Theorem 6.5.7]{Beardon 1991} 
or \cite[Appendix A.2]{Kawahira 2009}.\footnote{
{\it A priori} the error term is $O(w^{\nu+1})$, 
but here it is refined to be $O(w^{2\nu})$.}
The set of $w$'s with $\arg w^\nu=0$ (resp. $\arg w^\nu =\pi$) 
determines {\it the repelling} 
(resp. {\it attracting})  {\it directions} of this parabolic point. 
Note that the Fatou component $\Delta$ is invariant under 
$f_{c_1}^{k\nu}$,
and it contains a unique attracting direction.
In particular, the sequence $f_{c_1}^{k\nu m}(0)~(m \in \N)$ 
is contained in $\Delta$ and 
converges to $q_{c_1}$ as $m \to \infty$
tangentially to the attracting direction.

\begin{figure}[htbp]
\begin{center}
\includegraphics[width=.7\textwidth]{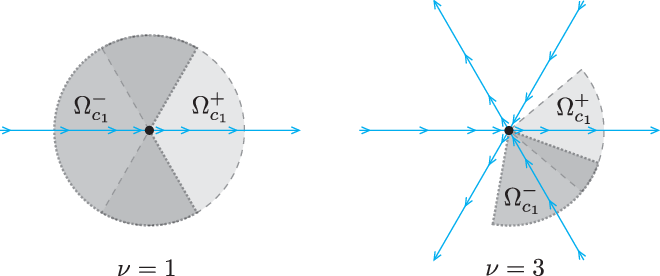}
\end{center}
\caption{\small We choose a pair of repelling and attracting petals.
Their intersection has two components when $\nu=1$.}
\label{fig_petals}
\end{figure}

Set 
\begin{align*}
\Omega_{c_1}^+&:=
\braces{
z=\psi_{c_1}^{-1}(w) 
\in \C \st -\frac{2\pi}{3\nu} \le \arg w \le \frac{2\pi}{3\nu},~0 < |w|<r
},\\
\Omega_{c_1}^-&:=
\braces{
z=\psi_{c_1}^{-1}(w) \in \C \st 
-\frac{5\pi}{3\nu} \le \arg w \le -\frac{\pi}{3\nu}
,~0 < |w|<r}
\end{align*}
for some sufficiently small $r>0$ such that 
$\Omega_{c_1}^+$ and $\Omega_{c_1}^-$ are a pair of repelling and attracting petals 
with $\Omega_{c_1}^+ \cap \Omega_{c_1}^- \neq \emptyset$.
(See Figure \ref{fig_petals}.)
By multiplying a $\nu$-th root of unity to the local coordinate $w=\psi_{c_1}(z)$ if necessary, 
we may assume that the attracting petal $\Omega_{c_1}^-$ is contained in $\Delta$.

For the coordinate $w=\psi_{c_1}(z)$, 
we consider an additional coordinate change 
$w \mapsto W=-1/(\nu w^\nu)$. In this $W$-coordinate,
the action of $f_{c_1}^{k\nu}$ on each petal is 
$$
W \mapsto W+1+ O(W^{-1}).
$$
By taking a smaller $r$ if necessary,
there exist 
conformal mappings 
$\phi_{c_1}^ + :\Omega_{c_1}^ +  \to \C$
and 
$\phi_{c_1}^-:\Omega_{c_1}^- \to \C$
such that $\phi_{c_1}^\pm (f_{c_1}^{k\nu}(z)) =\phi_{c_1}^\pm(z)+1$
which are unique up to adding constants.
(We will normalize them later.) 
We call $\phi_{c_1}^\pm$ the {\it Fatou coordinates}.

\parag{Perturbed Fatou coordinates.}
When $\nu=1$ (equivalently, $\mu_{c_1}=1$), 
the parabolic fixed point $q_{c_1}$ of $f_{c_1}^k$ 
splits into two distinct fixed points of 
$f_c^k$ for each $c \neq c_1$. 
To describe this bifurcation, 
it is convenient to use a parameter $u$ that satisfies $c=c_1+u^2$. 
It is known that there are two holomorphic functions 
$q_+(u)$ and $q_-(u)$ defined near $0$ such that 
$f_{c_1+u^2}^k(q_\pm(u))=q_\pm(u)$; 
$q_{c_1}=q_+(0)=q_-(0)$; and their multipliers satisfy
$$
\mu_\pm(u) := (f_{c_1+u^2}^k)'(q_\pm(u))=1 \pm A_0 u + O(u^2)
$$
for some $A_0 \neq 0$. 
(See \cite[Expos\'e XI]{DH Orsay}, \cite[Theorem 1.1 (c)]{Tan Lei 2000},
 or the primitive case of \cite[Lemma 4.2]{Milnor 2000}.)
Note that the maps $u \mapsto \mu_\pm(u)$ are univalent near $u=0$
and hence locally invertible. 

\begin{figure}[htbp]
\begin{center}
\includegraphics[width=.55\textwidth]{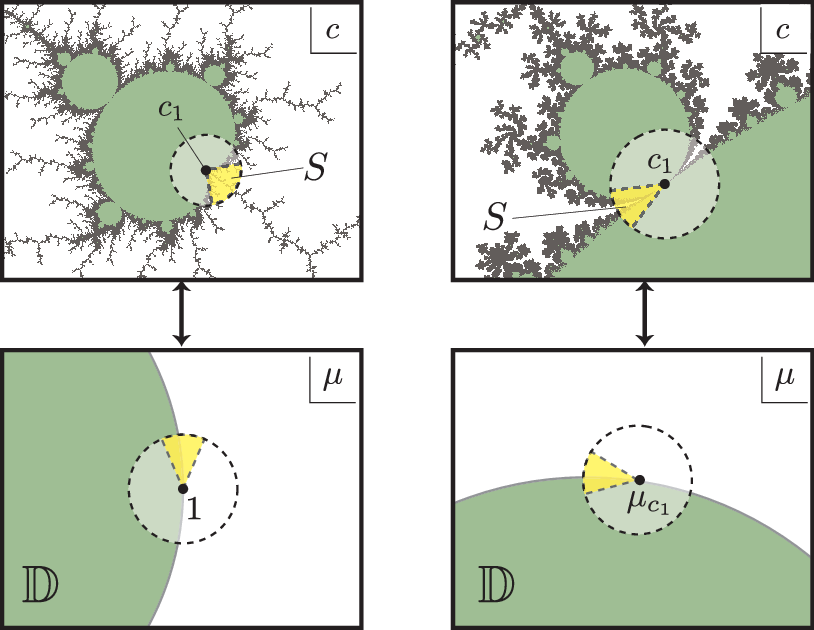}
\end{center}
\caption{\small The sector $S$ for $\nu=1$ (left) and $\nu \ge 2$ (right).}
\label{fig_sectors}
\end{figure}

For $r > 0$ we define 
a sector $S_\mu(r) \subset \C$ attached to $1$ by
$$
S_\mu(r):=\braces{\mu \in \C \st 
0<|\mu-1|<r
~~~
\text{and}~~~
\abs{\arg (\mu -1)-\frac{\pi}{2}}<\frac{\pi}{8}
}.
$$
We choose a sufficiently small $r_0>0$ such that 
the set
$$
S:=\braces{c=c_1+u^2  \st \mu_+(u) \in S_\mu(r_0)}
$$
is contained in $\Lambda$ 
and that the correspondence between $\mu=\mu_+(u) \in S_\mu(r_0)$ 
and $c=c_1+u^2 \in S$ is one-to-one.
See Figure \ref{fig_sectors} (left).
We may regard the parameter $u$ that 
mediates this one-to-one correspondence 
as a holomorphic branch of $\sqrt{c-c_1}$ over $S$.
We may also regard 
\begin{align}
q_c&:=q_+(u)=q_+(\sqrt{c-c_1})\qquad \text{and} \notag\\
\mu_c&:=\mu_+(u)=1+A_0 \sqrt{c-c_1} + O(c-c_1) \label{eq_A_0}
\end{align}
as a fixed point of $f_c^k$ and its multiplier that depend
holomorphically on $c \in S$.

When $\nu\ge 2$ (equivalently, $\mu_{c_1} \neq 1$),
the parabolic fixed point $q_{c_1}$ of $f_{c_1}^k$ 
splits into one fixed point $q_c$ and 
a cycle of period $\nu$ of $f_c^k$ for each $c \neq c_1$. 
By the implicit function theorem,
$q_c$ and its multiplier $\mu_c:=(f_c^k)'(q_c)$ depend holomorphically on $c$ near $c_1$,
and it is known that 
\begin{equation}\label{eq_B_0}
\mu_c := (f_c^k)'(q_c) = \mu_{c_1} \, \paren{1+B_0 (c-c_1)+ O((c-c_1)^2)}
\end{equation}
for some constant $B_0 \neq 0$.
(See \cite[Expos\'e XI]{DH Orsay}, \cite[Theorem A.1 (c)]{Tan Lei 2000},
or the satellite case of \cite[Lemma 4.2]{Milnor 2000}.)
Note that the map $c \mapsto \mu_c$ is univalent near $c_1$
and hence locally invertible. 

We choose a sufficiently small $r_0>0$ such that the set
$$
S:=\braces{c \in \C \st \frac{\mu_c}{\mu_{c_1}} \in S_\mu(r_0)}
$$
is contained in $\Lambda$ 
and that 
the correspondence between $\mu=\mu_c \in  \mu_{c_1} \times S_\mu(r_0)$ 
and $c \in S$ is one-to-one.
See Figure \ref{fig_sectors} (right).
We call $S$ a {\it sector} attached to $c_1 \in \partial M_{s_0}$.

Now we consider general $\nu$ ($=1$ or $\ge 2$).
By taking a sufficiently small $r_0$, 
we may assume that for each $c \in S$
there exists a holomorphic local coordinate $w=\psi_c(z)$ 
near $q_{c}$ with $\psi_c(q_{c})=0$ such that
$$
\psi_c \circ f_{c}^{k\nu} \circ \psi_c^{-1}(w)
=\mu_c^\nu w\,(1 + w^{\nu}+O(w^{2\nu})),
$$
where $\mu_c^\nu \to 1$ and $\psi_c \to \psi_{c_1}$ 
uniformly as $c \in S$ tends to $c_1$.
See \cite[Appendix A.2]{Kawahira 2009}. 
By a further $\nu$-fold coordinate change $W=-\mu_c^{\nu^2}/(\nu w^\nu)$,
the action of $f_c^{k \nu} $ is 
$$
W \mapsto \mu_c^{-\nu^2} \,W +1+O(W^{-1}),
$$ 
where $W=\infty$ is a fixed point 
with multiplier $\mu_c^{\nu^2}$ 
that corresponds to the fixed point $q_c$ of $f_c^{k\nu}$.
There is another fixed point of the form $W=1/(1-\mu_c^{-\nu^2})+O(1)$
with multiplier close to $\mu_c^{-\nu^2}$
on each branch of the $\nu$-fold coordinate.
Note that  
\begin{align*}
\mu_c^{\pm \nu^2}&=
1 \pm A_0\sqrt{c-c_1}+O(c-c_1)
\qquad \text{or}\\
\mu_c^{\pm \nu^2}&=1 \pm \nu^2 B_0(c-c_1)+O((c-c_1)^2)
\end{align*}
according to $\nu=1$ or $\nu \ge 2$ by (\ref{eq_A_0}) and (\ref{eq_B_0}).

\begin{figure}[htbp]
\begin{center}
\includegraphics[width=.44\textwidth]{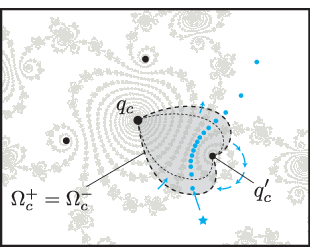}
\end{center}
\caption{\small A typical behavior of the critical orbit 
by $f_c^{k \nu }$ near $q_c$ for $\nu =3$.}
\label{fig_implosion}
\end{figure}

It is known that for each $c$ in $S \cup \{c_1\}$ 
(by taking a smaller $r_0$ if necessary),
there exist (perturbed) {\it Fatou coordinates} 
$\phi_c^ + :\Omega_c^ +  \to \C$ 
and
$\phi_c^-:\Omega_c^- \to \C$
satisfying the following conditions
(\cite{Lavaurs 1989}, \cite{DSZ 1997} and \cite[Proposition A.2.1]{Shishikura 1998}):

\begin{itemize}
\item
For any $c \in S$,
both $\partial \Omega_c^+$ and $\partial \Omega_c^-$ contain two fixed points $q_c$ and $q_c'$ 
of $f_c^{k\nu}$ that converge to $q_{c_1}$ as $c \in S$ tends to  $c_1$.
\item 
For any $c \in S \cup \{c_1\}$, 
$\phi_c^\pm$ is a conformal map from a domain $\Omega_c^\pm$ onto the image in $\C$ that satisfies
$\phi_c^\pm(f_c^{k\nu}(z))=\phi_c^\pm(z)+1$ 
if both $z$ and $f_c^{k\nu}(z)$ are contained in $\Omega_c^\pm$.
\item
(Holomorphic dependence) 
Every compact set $E$ in $\Omega_{c_1}^\pm$ 
is contained in $\Omega_c^\pm$ for $c \in S$ sufficiently close to $c_1$,
and $\phi_c^\pm(z)$ depends holomorphically on $c \in S$ near $c_1$ for each $z \in \Omega_{c_1}^\pm$.
\end{itemize}
Indeed, Fatou coordinates $\phi_c^\pm$ for each $c \in S \cup \{c_1\}$ are uniquely determined up to composition with a translation of $\C$. 
We will give a specific normalization for them 
by using the critical orbits in Step (P2) (Section 5). 
We can also arrange the domains $\Omega_c^\pm$
such that $\Omega_c^+=\Omega_c^-=:\Omega_c^\ast$ 
for each $c \in S$
(Figure \ref{fig_implosion}).
Hence for each $z \in \Omega_c^\ast$, 
$$
\tau(c):=\phi_c^+(z)-\phi_c^-(z) \in \C
$$
is defined and independent of $z$. 
The function $\tau:S \to \C$
is called the {\it lifted phase}.
Note that the value $\tau(c)$ is determined 
by the unique normalized Fatou coordinates
associated with the analytic germ $f_c^{k\nu}$,
and it does not depend on the choice of the parametrization.
\footnote{
In \cite{Douady 2000}, the lifted phase 
for $\nu=1$ is described 
in terms of the normalized germ
$f_\mu (z)=z+z^2+\mu+\cdots~(\mu \to 0)$.
In this case the multiplier for two fixed points of $f_\mu$ are 
$1 \pm 2 \sqrt{\mu}i\,(1+O(\mu))$ and $\tau=-\pi/\sqrt{\mu}+O(1)$
as $\mu \to 0$. 
In \cite{DSZ 1997}, they use $\al=(\nu/2\pi i)\log (\mu_c/\mu_{c_1}-1)$
(so that $\mu_c=\exp(2 \pi i (\nu'+\al)/\nu)$)
to parametrize the germs. 
Any parameterizations are analytically equivalent and they determine
the same value $\tau$ under appropriate normalizations.
}
It is known that if $\nu=1$, then
$$
\tau(c)= -\frac{2\pi i}{A_0\sqrt{c-c_1}}+O(1)
$$
as $c \in S$ tends to $c_1$, where $A_0$ is given in (\ref{eq_A_0}).
Similarly if $\nu \ge 2$, we obtain
$$
\tau(c)= -\frac{2 \pi i}{\nu^2 B_0(c-c_1)}+O(1)
$$
as $c \in S$ tends to $c_1$,
where $B_0$ is given in (\ref{eq_B_0}). 
In both cases it can be also shown that $\tau(c)$ is univalent on $S$ if 
$S$ is sufficiently small.

\section{Step (P1): Definitions of $U_c$, $U_c'$ and $V_c$}

Now we start the main steps of the proof of Theorem \ref{thm_A'}.
Let us briefly recall the notation: 
The parameter $s_0 \neq 0$ is superattracting
for which $P_{s_0}^p(0)=0$ with $p \ge 2$.
The small Mandelbrot set $M_{s_0}$ with center $s_0$,
and the family 
$
\{
f_c := P_c^p|_{\widetilde{U}_c'} : \widetilde{U}_c' \to \widetilde{U}_c
\}_{c \, \in \, \Lambda}
$
of quadratic-like maps are associated with it.
For a given parabolic parameter $c_0 \in \partial M$,
we have considered another parabolic parameter 
$c_1=s_0 \perp c_0 \in M_{s_0}$ and
its perturbation in a sector $S$.

\indent
For a technical reason, 
{\it 
we first assume that 
the parabolic parameter $c_1=s_0 \perp c_0 \in M_{s_0}$ 
belongs to $\Lambda$}.
Note that this assumption only excludes the case 
where $M_{s_0}$ is a satellite small Mandelbrot set and $c_0=1/4$. This case will be discussed separately.

Under this assumption, we shall construct a 
family $\bs{f}= \{f_c:U_c' \to U_c\}_{c\, \in \,S \cap \Lambda}$ 
of quadratic-like maps and a family  
$\bs{g}=\{g_c:V_c \to U_c\}_{c\, \in \,S\cap \Lambda}$ of isomorphisms.

We start with the ``global" dynamics of $P_{c_1}$:

\begin{lem}[cf. Lemma 4.1 of \cite{Kawahira-Kisaka 2018}]\label{lem_4.1'}
There exist Jordan domains $U_{c_1}$, $U_{c_1}'$ and $V_{c_1}$ with 
analytic boundaries and 
integers $N, \ j \ge 1$ which satisfy 
the following:
\begin{itemize}
\item[\rm (1)]
$0 \in U_{c_1}' \subset \widetilde{U}_{c_1}'$ and 
$f_{c_1}: U_{c_1}' \to U_{c_1}$ is a 
quadratic-like map.

\item[\rm (2)]
$g_{c_1} := P_{c_1}^N|_{V_{c_1}} : V_{c_1} \to U_{c_1}$ is an isomorphism
and $\overline{f_{c_1}^j(V_{c_1})} 
\subset U_{c_1} \smallsetminus \overline{U_{c_1}'}$.

\item[\rm (3)]
$\overline{V_{c_1}} \subset \Omega_{c_1}^+$.
\end{itemize}
There are infinitely many possible choices for $V_{c_1}$ 
(and the corresponding $N$ and $j$ associated with $V_{c_1}$), 
and $V_{c_1}$ can be chosen to have an arbitrarily small diameter. 
Also we can take $V_{c_1}$ arbitrarily
close to $q_{c_1} \in \partial \Omega_{c_1}^+$.
\label{def of U_{c_1} etc for parabolic case}
\end{lem}
\noin
Indeed, each choice of $V_{c_1}$ will determine a different 
copy of the decorated Mandelbrot set.

\parag{Proof.}
By shrinking $\widetilde{U}_{c_1}$ and $\widetilde{U}_{c_1}'$ slightly,
we can take Jordan domains $U_{c_1}$ and $U_{c_1}'$ 
with analytic boundaries which are neighborhoods of $J(f_{c_1})$ 
and $f_{c_1} : U_{c_1}' \to U_{c_1}$ is a quadratic-like map.

For $j =0,\,1,\, 2\, \ldots$ let
$$
  A_j := f_{c_1}^{-j}(U_{c_1} \smallsetminus \overline{U_{c_1}'}).
$$
Since $f_{c_1}$ is conjugate to $z^2$ on 
$U_{c_1}' \smallsetminus J(f_{c_1})$,
the annulus $A_j$ is uniformly close to $J(f_{c_1})$ and  
thus $A_j \cap \Omega_{c_1}^+ \ne \emptyset$ for every sufficiently large $j$.
Also since the ``global" Julia set $J(P_{c_1})$ of $P_{c_1}$ is a connected set containing the ``small" Julia set $J(f_{c_1})$, 
the annulus $A_0=U_{c_1} \smallsetminus \overline{U_{c_1}'}$
intersects with $J(P_{c_1})$ 
and so does $A_j$ for $j \ge 1$.
In particular, for any $\zeta_0 \in J(P_{c_1}) \cap A_0$, 
the set $f_{c_1}^{-j}(\{\zeta_0\})$ 
approximates $J(f_{c_1})$ as $j \to \infty$
in the Hausdorff topology.
Hence $A_j \cap \Omega_{c_1}^+$ contains 
a point $z_0 \in J(P_{c_1})$ 
arbitrarily close to the parabolic periodic point 
$q_{c_1} \in J(f_{c_1})$
for every sufficiently large $j$.
Let $B$ be any closed disk in $A_j \cap \Omega_{c_1}^+$
centered at this $z_0$ with a sufficiently small radius.

Note that the postcritical set of 
the map $P_{c_1}$ in $\mathbb C$ 
is contained in 
$U_{c_1}' \cup P_{c_1}(U_{c_1}') \cup \cdots \cup P_{c_1}^{p-1}(U_{c_1}')$.
There are two disjoint connected components 
$X:=P_{c_1}^{p-1}(U_{c_1}')$
and 
$-X:=\{-x \in \mathbb{C} ~:~ x \in X\}$
of $P_{c_1}^{-1}(U_{c_1})$, 
where $-X$ does not intersect with neither 
the postcritical set of $P_{c_1}$ nor the critical point $0$.
Hence for any $n \in \N$ and any connected component $V$ of
$P_{c_1}^{-n}(-X)$,
$P_{c_1}^{n+1}:V \to U_{c_1}$ is an isomorphism.

Since the inverse images of $-X$ in the dynamics of $P_{c_1}$
accumulate on any point in the Julia set $J(P_{c_1})$ of $P_{c_1}$ 
(by Montel's theorem), 
the shrinking lemma (\cite[p.86]{Lyubich-Minsky 1997} or 
\cite[Lem.2.9]{Cui-Tan 2018})
implies that we can find a component $V_{c_1}$ of 
$P_{c_1}^{-N+1}(-X)$ contained in the closed disk $B$ 
for some $N \in \N$. 
This gives a desired isomorphism $g_{c_1}:=P_{c_1}^N:V_{c_1} \to U_{c_1}$.\\
\hfill
\QED ~{\small (Lemma \ref{lem_4.1'})}
\\

\parag{Definition of the families
$\bs{f}$ and $\bs{g}$.}
For $ \vep'>0$ given in the statement of Theorem \ref{thm_A'},
we may assume that the sector $S$ attached to $c_1$
(defined in the previous section)
is contained in $D(c_1, \vep') \cap \Lambda$
by taking a smaller $r_0>0$ in the definition of $S$.

\begin{figure}[htbp]
\begin{center}
\includegraphics[width=.60\textwidth]{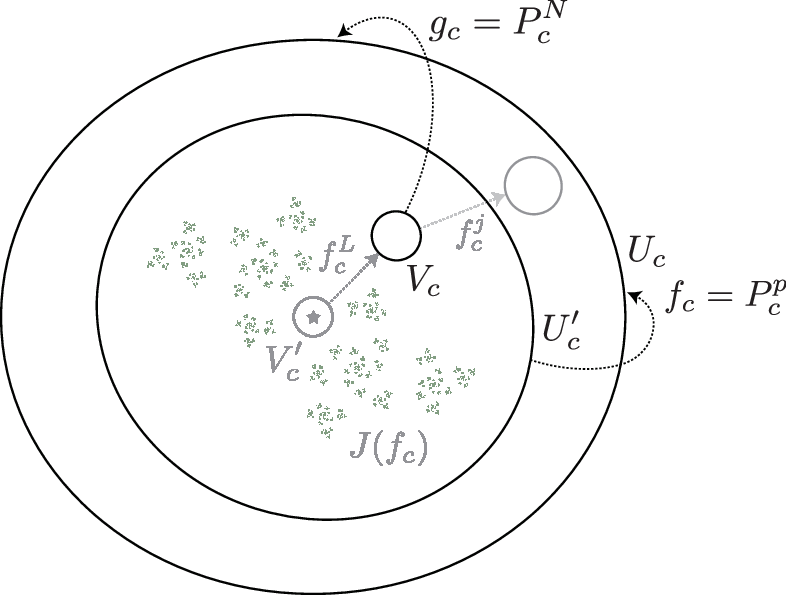}
\caption{\small
\small
The Jordan domains $U_c$, $U_c'$, $V_c$, and $V'_c$ for parameters $c \in S \smallsetminus \{c_1\}$. When $c = c_1$, the domains $U_{c_1}$, $U_{c_1}'$, and $V_{c_1}$ are arranged as in the figure, but the small Julia set $J(f_{c_1})$ is connected. The additional domain $V_c'$ will be defined later and exists only for $c \in S \smallsetminus \{c_1\}$.
}\label{UV}
\end{center}
\end{figure}

Let $U_c := U_{c_1}$ for each $c \in S \subset \Lambda$.
By taking a smaller $S$ if necessary,
we obtain a quadratic-like map $f_c:U_c' \to U_c$
with $U_c':=f_c^{-1}(U_c) \subset \widetilde{U}_{c}'$,
a component $V_c$ of $P_c^{-N}(U_c)$ 
 close to $V_{c_1}$ satisfying 
 $\overline{f_c^j(V_c)} \subset U_c \smallsetminus \overline{U_c'}$,
and an isomorphism $g_c : V_c \to U_c \,(\equiv U_{c_1})$ 
for each $c \in S$.
See 
Figure \ref{UV}. 
We also define the {\it pre-critical point} $b_c$
by $b_c:=g_c^{-1}(0) \in V_c$. 
Hence we have done the construction of the families
$$
\bs{f}= \braces{f_c:U_c' \to U_c}_{c \, \in\, S}
\quad\text{and}\quad
\bs{g}=\braces{g_c:V_c \to U_c}_{c \, \in\, S}.
$$

\parag{Satellite roots.}
Now we deal with the remaining case 
where $M_{s_0}$ is a satellite small Mandelbrot set 
with renormalization period $p$ and $c_0 =1/4$.
(We say $c_1$ for this case is a {\it satellite root}.) 
Let $\Lambda$ be the simply connected domain associated with $M_{s_0}$. 
Since $M_{s_0}$ is satellite, 
the quadratic-like family 
$
\{ 
f_c = P_c^p|_{\widetilde{U}_c'} : \widetilde{U}_c' \to \widetilde{U}_c
\}_{c \, \in \,\Lambda}
$
excludes the parameter $c_1 =s_0 \perp (1/4) \notin \Lambda$.
However, by slightly modifying the notion of quadratic-like map 
at this parameter,
one can establish a version of Lemma 4.1 as follows.

Let $q_{c_1}$ be the fixed point of $P_{c_1}^{p}$ with multiplier $1$, 
whose petal number is $\nu \ge 2$. 
(Then the period of $q_{c_1}$ in the dynamics of $P_{c_1}$ 
is $p'=p/\nu$.)
Let $\Omega_{c_1}^+$ be the repelling petal attached to $q_{c_1}$ as in Figure \ref{fig_petals} (right). 
Then we have:

\parag{Lemma \ref{lem_4.1'} for the Satellite Roots.} 
{\it 
There exist Jordan domains $U_{c_1}$, $U_{c_1}'$ and $V_{c_1}$ with 
analytic boundaries and 
integers $N, \ j \ge 1$ which satisfy 
the following:
\begin{itemize}
\item[\rm (1)]
$0 \in {U}_{c_1}' \subset {U}_{c_1}$,
$\partial {U}_{c_1}' \cap \partial {U}_{c_1}=\braces{q_{c_1}}$,
and $f_{c_1}:=P_{c_1}^p|_{{U}_{c_1}'}:{U}_{c_1}' \to {U}_{c_1}$
is a proper branched covering of degree two.
\item[\rm (2)]
$g_{c_1} := P_{c_1}^N|_{V_{c_1}} : V_{c_1} \to U_{c_1}$ is an isomorphism
and $\overline{f_{c_1}^j(V_{c_1})} 
\subset U_{c_1} \smallsetminus \overline{U_{c_1}'}$.
\item[\rm (3)]
$\overline{V_{c_1}} \subset \Omega_{c_1}^+$. 
\end{itemize}
There are infinitely many possible choices for $V_{c_1}$ 
(and the corresponding $N$ and $j$ associated with $V_{c_1}$), 
and $V_{c_1}$ can be chosen to have an arbitrarily small diameter. 
Also we can take $V_{c_1}$ arbitrarily
close to $q_{c_1} \in \partial \Omega_{c_1}^+$.
}

\medskip
\noindent
The proof is analogous to those of Lemmas 4.1. 
However, to obtain (1) we need 
an important idea 
of the proof of the Douady-Hubbard tuning theorem 
\cite[\S 3]{Haissinsky 2000} 
which is essential in the construction of $\Lambda$.
Here we only give a sketch:
For each $c \in \Lambda$ there exists a repelling fixed point $q_c$ of $P_{c}^{p}$ that depends holomorphically on $c$ and 
$q_c \to q_{c_1}$ as $c \in \Lambda$ tends to $c_1$.
We define a Jordan domain $U_c$ by adding 
a small disk centered at $q_c$ to $U_{c_1}$.
(We also modify $U_c$ slightly such that 
$\partial U_c$ is an analytic Jordan curve 
that moves holomorphically with respect to $c$.)
Then for any $c \in \Lambda$ sufficiently close to $c_1$,
we have a quadratic-like map $f_c:U_c' \to U_c$
and an isomorphism $g_c:V_c \to U_c$
where $U_c'$ is a connected component of $P_c^{-p}(U_c)$ 
with $\overline{U_c'} \subset U_c$
and $V_c$ is a connected component of $P_c^{-N}(U_c)$
that is close to $V_{c_1}$ and
$\overline{f_{c}^j(V_{c})} \subset U_{c} \smallsetminus \overline{U_{c}'}$. 
Since $q_{c_1}$ is parabolic with $\nu \ge 2$ petals,
we take the sector $S$ attached to $c_1$ as in Figure \ref{fig_sectors} (right).
By taking $S$ with sufficiently small radius, we obtain
the families $\{f_c:U_c' \to U_c\}_{c \,\in \,S \,\cap \,\Lambda}$ 
and $\{g_c:V_c \to U_c\}_{c \,\in \,S \,\cap\, \Lambda}$
over the set $S \cap \Lambda$ together with Fatou coordinates 
and lifted phase. 
In conclusion, we have constructed the families
$$
\bs{f}= \braces{f_c:U_c' \to U_c}_{c \, \in\, S \cap \Lambda}
\quad\text{and}\quad
\bs{g}= \braces{g_c:V_c \to U_c}_{c \, \in\,  S \cap \Lambda}.
$$

\parag{Remark.}
We have $S \not\subset \Lambda$ for the satellite root $c_1$,
but $S \subset \Lambda$ for the other cases. 
Hence we regard $\{f_c\}$ and $\{g_c\}$ as families defined over $S \cap \Lambda$ for all cases.

\medskip

\section{
Step (P2): Construction of the quadratic-like family 
$\boldsymbol{G}$}

We shall construct the second quadratic-like family 
$$
\boldsymbol{G} 
:= 
\{ G_c : V_c' \to U_c 
   \}_{c \,\in\, W} 
$$
such that $V_c' \subset U_c'$ and $W \subset S \cap \Lambda$.

\parag{\bf Normalization of the Fatou coordinates.}
Recall that we have (perturbed) Fatou coordinates 
$\phi_c:\Omega_c^\pm \to \C$ for each $c \in S \cup \{c_1\}$
as given in Section 3. 
By taking a smaller $S$ if necessary, 
we may normalize them such that:
\begin{itemize}
\item
$\overline{V_c} \subset \Omega_c^+$ for any $c \in S \cup \{c_1\}$. 
\item
There exists an $m \in \N$
such that $f_{c}^{k\nu m}(0) \in \Omega_c^-$ 
for any $c \in S \cup \{c_1\}$,
and $\phi_c^-$ is normalized such that $\phi_c^-(f_{c}^{k\nu m}(0))=m$.
\item
$\phi_c^+(b_c)=0$, where $b_c=g_c^{-1}(0) \in V_c$ is the pre-critical point.
\end{itemize}
Recall also that
we may arrange the domains $\Omega_c^\pm$
such that $\Omega_c^+=\Omega_c^-=:\Omega_c^\ast$ 
for each $c \in S$.

\parag{\bf Definition of $W$.}
Now for each $n \ge m$, define
$$
W = W_n 
:= 
\{ c \in S \ | \ f_c^{k\nu i}(0) \in \Omega_c^* \ \text{ for } \ i=m, \dots, n-1
\ \text{ and } \ 
  f_c^{k\nu n}(0) \in V_c \},
$$
that is, we consider the parameter $c$ such that the orbit of $0$
by $f_c^{k\nu}$ hits $V_c$.

\medskip

\begin{lem}[cf. Lemma 4.2 of \cite{Kawahira-Kisaka 2018}]
\label{lem_4.2}
By shrinking $U_c \equiv U_{c_1}$ slightly, the set $W=W_n$ is a non-empty
Jordan domain with analytic boundary for every 
sufficiently large $n$. 
Moreover, there 
exists an $s_n \in W_n$ such that $f_{s_n}^{k\nu n}(0) = b_{s_n}$,
which implies $g_{s_n} \circ f_{s_n}^{k\nu n}(0) = P_{s_n}^{k\nu np+N}(0) = 0$ 
and hence $P_{s_n}$ has a superattracting periodic point. 
\end{lem}

\medskip

\parag{Proof.}
We observe the dynamics of $f_c^{k\nu}$ near $q_c$ through
the perturbed Fatou coordinate 
$\phi_c^+ : \Omega_c^+=\Omega_c^\ast \to \C$ of $q_c$. 
Let
$$
  \wt{V}_c := \phi_c^+(V_c),
$$
then $c \in W_n$ if and only if $\phi_c^+(f_c^{k\nu n}(0)) \in \wt{V}_c$.
By the normalization of $\phi_c^+$, we have
$$
  \tau(c) = \phi_c^+(f_c^{k\nu n}(0)) - \phi_c^-(f_c^{k\nu n}(0))
= \phi_c^+(f_c^{k\nu n}(0)) - n.
$$
Hence it follows that $c \in W_n$ if and only if 
$$
  \tau(c) + n \in \wt{V}_c.
$$
Next take a Riemann map
$$
  u : U_c \equiv U_{c_1} \to \D, \quad u(0) = 0
$$
and define
$$
  v(c, \zeta) := \phi_c^+ \circ (u \circ g_c)^{-1}(\zeta), \quad 
  \zeta \in \D \ \text{ with } \ v(c, 0)=0, 
$$
which is the inverse of a Riemann map $u \circ g_c \circ (\phi_c^+)^{-1}$
of $\wt{V}_c$. 
Now we solve the equation with respect to the variable $c$
\begin{equation}
  \tau(c)+n = v(c,\zeta)
\label{eqn for v(c,zeta) for parabolic case}
\end{equation}
for each fixed $\zeta \in \D$. 
More precisely, we show that 
there exists an $n_0 \in \N$ such that for every $n \geq n_0$ and
$\zeta \in \D$ the equation (\ref{eqn for v(c,zeta) for parabolic case}) has a unique
solution.

\noin
{\bf Case 1 : $\nu = 1$. }
Since
$$
\tau(c)
=
-\frac{2\pi i}{A_0\sqrt{c - c_1}} + h(c), \quad h(c)=O(1) \quad (c \to c_1),
$$
the equation (\ref{eqn for v(c,zeta) for parabolic case}) 
can be rewritten as
\begin{equation}
F(c,\zeta) + G(c,\zeta) = 0,
\label{eqn F+G=0 for parabolic case}
\end{equation}
where
$$
  F(c,\zeta) :=  -\frac{2\pi i}{A_0\sqrt{c - c_1}} + n - v(c_1,\zeta), \quad
  G(c,\zeta) := h(c) - \big( v(c,\zeta)- v(c_1,\zeta) \big).
$$
The equation $F(c,\zeta)=0$ has a unique solution
$$
  c = c_n(\zeta) := c_1 -\frac{4\pi^2}{A_0^2(n - v(c_1,\zeta))^2}.
$$
Let
$$
  r_n(\zeta) := \bigg| -\frac{4\pi^2}{A_0^2(n - v(c_1,\zeta))^2} \bigg| 
= O(n^{-2})
$$
and take any $\beta$ with $0 < \beta < 1/2$. Consider 
(\ref{eqn F+G=0 for parabolic case}) in the disk 
$D(c_n(\zeta), r_n(\zeta)^{1+\beta})$. Since it is easy to see that
$$
  |F(c,\zeta)| = O(r_n(\zeta)^{\beta-1/2}) = O(n^{1-2\beta}), \quad 
  |G(c,\zeta)| = O(1)
$$
on the boundary $C$ of this disk, 
we have $|F(c,\zeta)| > |G(c,\zeta)|$ on $C$ for sufficiently large $n$. 
By Rouch\'e's theorem (\ref{eqn F+G=0 for parabolic case}) has
a unique solution $c = \check{c}_n(\zeta)$ in $D(c_n(\zeta), r_n(\zeta)^{1+\beta})$.
By using this solution, we can write
$$
  W_n = \{ \check{c}_n(\zeta) \in \C \ | \ \zeta \in \D \}.
$$

\medskip

\parag{Claim.} \ 
{\it 
{\rm (1)} $\check{c}_n : \D \to W_n$ is holomorphic.

\noindent
{\rm (2)} For every $r \in (0,1)$, $\check{c}_n$ is univalent on $\D(r)$ for every
sufficiently large $n$.
}

\medskip

\parag{Proof.} 
(1) By the argument principle (\cite[p.153, (49)]{Ahlfors 1978}), for each $\zeta \in \D$ we have
$$
  \check{c}_n(\zeta)
= \frac{1}{2 \pi i}
  \int_C  H(c,\zeta) c \cdot dc, 
$$
where
$$
H(c,\zeta) 
:= 
\frac{\frac{\partial}{\partial c} \big(F(c,\zeta) + G(c,\zeta)\big)}
{F(c,\zeta)+G(c,\zeta)}, \quad
C = \{ z \ | \ |c-c_n(\zeta)|=r_n(\zeta)^{1+\beta} \}
$$
Hence if $|\Delta \zeta| \ll 1$ and $c_n(\zeta+\Delta \zeta) \in \text{int}(C)$,
we have
$$
  \check{c}_n(\zeta+\Delta\zeta)
= \frac{1}{2 \pi i}
  \int_C  H(c,\zeta+\Delta\zeta) c \cdot dc.
$$
Then it follows that $H$ is holomorphic with respect to $\zeta$ and hence
$\check{c}_n(\zeta)$ is holomorphic in a neighborhood of $\zeta$. Thus
$\check{c}_n : \D \to W_n$ is holomorphic.

\noindent
(2) Since
$$
  \tau(\check{c}_n(\zeta)) + n = v(\check{c}_n(\zeta), \zeta)
$$
and $\tau(c)$ is univalent, it is enough to show that
$v(\check{c}_n(\zeta), \zeta)$ is univalent on $\overline{\D(r)}$ for
sufficiently large $n$. Note that 
$v(\check{c}_n(\zeta), \zeta) \to v(c_1, \zeta) \ (n \to \infty)$
uniformly on $\overline{\D(r)}$.

Set $v_n(\zeta):=v(\check{c}_n(\zeta), \zeta)$
and $v(\zeta):=v(c_1, \zeta)$ for brevity.
Now suppose on the contrary that
$$
  v_n(\zeta_n) = v_n(\zeta_n')
$$
for some $\zeta_n, \ \zeta_n' \in \overline{\D(r)}, \ \zeta_n \ne \zeta_n'$,
where $n$ ranges over a subsequence $\{ n_k \}_{k=1}^\infty$. By taking a
further subsequence, we may assume that
$$
  \zeta_n \to \hat{\zeta}, \quad \zeta_n' \to \hat{\zeta}' \quad \text{ for }
  \ n=n_k, \ k \to \infty.
$$
Then by passing to the limit in $v_n(\zeta_n) = v_n(\zeta_n') \ (n=n_k)$, 
we have $v(\hat{\zeta}) = v(\hat{\zeta}')$. By the univalence of $v$
it follows that $\hat{\zeta} = \hat{\zeta}'$.

Let $C_0 := v'(\hat{\zeta}) \ne 0$, then there exists a $\delta > 0$
such that
$$
  |v'(\zeta)-C_0| \leq \frac{|C_0|}{4} 
\quad \text{on} \ \D(\hat{\zeta}, \delta).
$$
Since $v_n \to v$ uniformly, we have
$$
  |v_n'(\zeta)-C_0| \leq \frac{|C_0|}{2} 
\quad \text{on} \ \D(\hat{\zeta}, \delta) \ \text{for} \ n \gg 0.
$$
Hence for $\zeta, \ \zeta' \in \D(\hat{\zeta}, \delta)$ we have
\begin{eqnarray*}
\big| \{ v_n(\zeta) - v_n(\zeta') \} - C_0(\zeta-\zeta') \big|
& = &
\Bigg| \int_\zeta^{\zeta'} (v_n'(\zeta) - C_0) d\zeta \Bigg| \\
& \leq &
\frac{|C_0|}{2} |\zeta-\zeta'|
\end{eqnarray*}
It follows that
$$
   \frac{|C_0|}{2} |\zeta-\zeta'|
\leq |v_n(\zeta) - v_n(\zeta')| 
\leq 
\frac 32 |C_0| |\zeta-\zeta'|.
$$
In particular, $v_n$ is injective on $\D(\hat{\zeta}, \delta)$ for
$n \gg 0$. However, $\zeta_n, \ \zeta_n' \in \D(\hat{\zeta}, \delta)$ 
for $n \gg 0$ and this is a contradiction.
\QED (Claim)

\medskip
By shrinking $U_c \equiv U_{c_1}$ slightly and using the Riemann map $u$ 
of the original $U_c \equiv U_{c_1}$, the boundary of the new 
$U_c \equiv U_{c_1}$ is parametrized as $u^{-1}(\gamma(t))$, 
where 
$\gamma(t) = re^{2\pi it} \in \D$ with $t \in [0,1]$ 
and $r \in (0,1)$ sufficiently close to $1$. 
Then 
$\partial \wt{V}_c$ is parameterized as $v(c, \gamma(t))$ and hence 
$\partial W_n$ (for the new $W_n$) is parameterized as 
$\check{c}_n(\gamma(t))$ by using the solution $\check{c}_n(\zeta)$ for 
the equation (\ref{eqn for v(c,zeta) for parabolic case}). Clearly this is an analytic Jordan curve
and $W_n$ is the image of $\D(r)$ by $\check{c}_n(\zeta)$. This shows that 
$W_n$ is a non-empty Jordan domain with analytic boundary. 

In particular, let $s_n := \check{c}_n(0)$ then this satisfies
$\tau(s_n) + n = 0$. This means that
$f_{s_n}^{k\nu n}(0) = b_{s_n}$,
which implies $g_{s_n} \circ f_{s_n}^{k\nu n}(0) = P_{s_n}^{k\nu np + N}(0)= 0$. 
Hence $P_{s_n}$ has a superattracting periodic point. 
This completes the proof for Case 1.

\medskip

\noin
{\bf Case 2 : $\nu \geq 2$. }
The argument is completely parallel to the Case 1. In this case
the functions $\tau(c), \ F(c,\zeta), \ G(c,\zeta), \ c_n(\zeta)$ and $r_n(\zeta)$
are replaced with
\begin{eqnarray*}
& & \tau(c)
= 
-\frac{2 \pi i}{\nu^2 B_0(c-c_1)} + h(c), \quad h(c) = O(1) \quad (c \to c_1), \\
& & F(c,\zeta) 
= 
-\frac{2 \pi i}{\nu^2 B_0(c-c_1)} + n - v(c_1,\zeta), \quad
G(c,\zeta)
=  h(c) - \big( v(c,\zeta)- v(c_1,\zeta) \big), \\
& & c_n(\zeta) 
=
c_1 + \frac{2 \pi i}{\nu^2 B_0(n-v(c_1,\zeta))} \quad \text{and} \quad
r_n(\zeta) = \bigg| \frac{2\pi i}{\nu^2 B_0(n - v(c_1,\zeta))} \bigg| = O(n^{-1}).
\end{eqnarray*}
Then we have the estimates
$$
  |F(c,\zeta)| = O(r_n(\zeta)^{-1+\beta}) = O(n^{1-\beta}), \quad
  |G(c,\zeta)| = O(1)
$$
on $\partial D(c_n(\zeta), r_n(\zeta)^{1+\beta})$ for $0 < \beta < 1/2$.
The rest of the argument is the same as in Case 1 and hence 
the same conclusion follows also in Case 2. 
This completes the proof. 

\hfill
\QED ~{\small (Lemma \ref{lem_4.2})}

\parag{Definition of the family $\bs{G}$.} 
Now we take a sufficiently large $n$
such that the previous lemma holds.  
We call $s_n \in W_n$ the {\it center} of $W_n$.
Let 
$L = L_n  := k\nu n$ and
$V_c'=V'_{c,n}$ be the component of $f_c^{-L}(V_c)$ containing 
$0$ and define
$$
  G_c = G_{c,n} := g_c \circ f_c^L : V_c' \to U_c
$$
for $c \in W_n=\{c \in S \st f_c^L(0) \in V_c\}$.
The map
$f_c^L : V_c' \to V_c$ is a branched covering of degree 2 and
$g_c : V_c \to U_c$ is a holomorphic isomorphism. 
Since $V_c$ satisfies 
$\overline{f_c^j(V_c)} \subset U_c \sminus \overline{U_c'}$,
we have 
$V_c' \subset f_c^{-L}(V_c) \subset f_c^{-L-j}(U_c) \Subset U_c$.
Hence 
$G_c := g_c \circ f_c^L : V_c' \to U_c $ is a quadratic-like map
for each $c \in W_n$ and 
we define a quadratic-like family $\boldsymbol{G}$ by
$$
\boldsymbol{G} = \boldsymbol{G}_n 
 := 
\{ G_c : V_c' \to U_c \}_{c \,\in\, W_n}.
$$
See Figure \ref{UV}.

\section{
Step (P3): Proof for $\boldsymbol{G}$ 
being a Mandelbrot-like family} 

In this step we follow Douady and Hubbard's formulation
to obtain a (quasiconformal) copy of $M$
in the parameter space of a given quadratic-like family. 

\parag{Mandelbrot-like families.}
A family of holomorphic maps $\boldsymbol{h} = \{ h_\lambda \}_{\lambda \,\in\, W}$ 
is called a {\it Mandelbrot-like family} if the 
following (1)--(8) hold: 

\begin{itemize}
\setlength{\itemsep}{5pt}
\item[(1)]
$W \subset \C$ is a Jordan domain with $C^1$ boundary $\partial W$.

\item[(2)]
There exists a family of maps 
$\Theta = \{ \Theta_\lambda \}_{\lambda \in W}$ such that
for every $\lambda \in W$, $\Theta_\lambda : \overline{A(R, R^2)} \to \C$ is a 
quasiconformal embedding for some $R > 1$, which is independent of $\lambda$,
and that $\Theta_\lambda(Z)$ is 
holomorphic in $\lambda$ for every $Z \in \overline{A(R, R^2)}$.

\item[(3)]
Define $C_\lambda := \Theta_\lambda(\partial D(R^2)), \
C'_\lambda := \Theta_\lambda(\partial D(R))$
and let 
$U_\lambda \ (\text{resp.} \ U'_\lambda)$ be the Jordan domain bounded by 
$C_\lambda \ (\text{resp.} \ C'_\lambda)$. Then
$h_\lambda$ is holomorphic in a neighborhood of $\overline{U'_\lambda}$
and
$h_\lambda : U'_\lambda \to U_\lambda$ is a quadratic-like map with
a critical point $\omega_\lambda$ depending holomorphically on $\lambda$. 
Also let
$$
  {\mathcal U} := \{ (\lambda, z) \ | \ \lambda \in W, \ z \in U_\lambda \}, \quad
  {\mathcal U'} := \{ (\lambda, z) \ | \ \lambda \in W, \ z \in U'_\lambda \}
$$
then
$\boldsymbol{h} : {\mathcal U'} \to {\mathcal U}, \ 
(\lambda, z) \mapsto (\lambda, h_\lambda(z))$ is analytic and proper.

\item[(4)]
$\Theta$ is equivariant on the boundary, i.e.,
$\Theta_\lambda(Z^2) = h_\lambda(\Theta_\lambda(Z))$ for $Z \in \partial D(R)$.

\vskip 2mm

\noin
The family of maps 
$\Theta = \{ \Theta_\lambda \}_{\lambda \in W}$ satisfying the above conditions
(1)--(4) is called a {\it tubing}. 

\vskip 2mm

\item[(5)]
$\boldsymbol{h}$ extends continuously to a 
map from $\overline{{\mathcal U}'}$ to $\overline{\mathcal U}$ and 
the map
$(\lambda ,z) \mapsto (\lambda, \Theta_\lambda(Z))$ extends 
continuously to a map from
$\overline{W} \times \overline{A(R, R^2)}$ to $\overline{\mathcal U}$
such that $\Theta_\lambda$ is injective on $A(R, R^2)$ for $\lambda \in \partial W$.

\item[(6)]
The map $\lambda \mapsto \omega_\lambda$ extends continuously to $\overline{W}$.

\item[(7)]
$h_\lambda(\omega_\lambda) \in C_\lambda$ for $\lambda \in \partial W$.

\item[(8)]
{\it The one turn condition}: 
when $\lambda$ ranges over $\partial W$ making one turn, then the vector
$h_\lambda(\omega_\lambda) - \omega_\lambda$ makes one turn around 0.
\end{itemize}
Now let $M_{\boldsymbol{h}}$ be the {\it connectedness locus} of the family
$\boldsymbol{h} = \{ h_\lambda \}_{\lambda \,\in\, W}$:
$$
M_{\boldsymbol{h}} := \{ \lambda \in W \ | \ K(h_\lambda)  \text{ is connected} \}
= \{ \lambda \in W \ | \ \omega_\lambda \in K(h_\lambda) \}.
$$
Douady and Hubbard (\cite[Chapter IV]{Douady-Hubbard 1985}) showed that there
exists a homeomorphism 
$$
  \chi : M_{\boldsymbol{h}} \to M.
$$ 
This is just a correspondence by the Douady-Hubbard straightening theorem, that is, for every
$\lambda \in M_{\boldsymbol{h}}$ there exist a unique
$c = \chi(\lambda) \in M$ such that $h_\lambda$ is hybrid equivalent	
to $P_c(z) = z^2 + c$. Furthermore they
showed that this
$\chi$ can be extended to a homeomorphism 
$\chi_\Theta : W \to W_M$ by using
$\Theta = \{ \Theta_\lambda \}_{\lambda \,\in\, W}$, 
where $W_M$ is a neighborhood of $M$ given by
$$
  W_M := \{ c \in \C \ | \ {\mathcal G}_M(c) < 2 \log R \}, 
  \quad {{\mathcal G}_M := }  \text{ the Green function of } M.
$$
In particular,
$\chi_\Theta(\lam)$ for  
$\lam \in W \sminus M_{\boldsymbol{h}}$ 
is defined in such a way that
\begin{equation}\label{eq_chi_Theta}
\Theta_\lam^{-1} 
\paren{h_\lam^k(\omega_\lam)}
=\braces{\Phi_M(\chi_\Theta(\lam)}^{2^{k-1}}
\end{equation}
for the unique $k \in \N$
with $h_\lam^k(\omega_\lam) \in U_\lam-U_\lam'$.
Also Lyubich showed that $\chi_\Theta$ is 
quasiconformal on any $W'$ with $W' \Subset W$ 
(\cite[p.366, THEOREM 5.5 (The QC Theorem)]{Lyubich 1999}).

\parag{Mandelbrot-like family with a ``decorated" tubing.}
Now recall that 
there exists a canonical homeomorphism
$\chi_{s_0}:\Lambda \to \chi_{s_0}(\Lambda)$
such that for any $c \in \Lambda$, 
$f_c:\tilde{U}_c' \to \tilde{U}_c$ 
is hybrid equivalent to some 
$P_{\alpha}$ with $\alpha=\chi_{s_0}(c)$
by the Douady-Hubbard tuning theorem. 
Here we will check:

\begin{lem}
\label{lem_G_is_M_like}
Let $W=W_n$, $s=s_n$, and $\sigma=\chi_{s_0}(s)$.
Then for sufficiently large $n \in \N$,
there exists an $R>1$ such that
the family 
$\bs{G}=
\braces{G_c : V_c' \to U_c}_{c \, \in\, W}$
is a Mandelbrot-like family
with a tubing 
$$
\Theta=
\{\Theta_c: \Bar{A(R,R^2)} \to \Bar{U_c} \sminus V_c'\}
_{c \, \in\, W}
$$
satisfying
$\Theta_c(\Gamma(\sigma))=J(f_c)$ for any $c \in W$.
\end{lem}

\parag{Proof.}
Suppose that $n$ is sufficiently large 
and the quadratic-like family $\boldsymbol{G}=\boldsymbol{G}_n$
over $W=W_n \subset S \cap \Lambda$ is defined as in Step (P2).
We construct a tubing 
$\Theta = \Theta_n = \{ \Theta_c \}_{c \,\in\, W_n}$ for 
$\boldsymbol{G}_n$
as follows: 
For $s_n \in W_n$,  since $f_{s_n}^L(0) \in V_{s_n}$ and 
$f_{s_n}^j(V_{s_n}) \subset U_{s_n} \smallsetminus \overline{U_{s_n}'}$, 
from Lemma \ref{def of U_{c_1} etc for parabolic case}, we have
$f_{s_n}^{L+j}(0) \notin U_{s_n}'$. 
It follows that $J(f_{s_n})$ is a Cantor set,
 which is quasiconformally homeomorphic to a quadratic Cantor
Julia set $J(P_{c_0+\eta_n})$ for some $\eta=\eta_n$ with 
$\chi_{s_0}(s_n)=c_0+\eta_n \notin M$ by the Douady-Hubbard straightening theorem.
Let $\Psi_{s_n}$ be the quasiconformal straightening map
that conjugates $f_{s_n}$ and $P_{c_0+\eta_n}$ 
defined on a neighborhood of $K(f_{s_n})$.
Then the image of $J(f_{s_n})$ by $\Psi_{s_n}$ is $J(P_{c_0+\eta_n})$.
Take an $R > 1$  such
that $J(P_{c_0+\eta_n}) \subset A(R^{-1/2}, R^{1/2})$. 
Let $\Gamma$ be the rescaled Julia set, that is, 
$$
  \Gamma := \Gamma_0(c_0+\eta_n) 
= J(P_{c_0+\eta_n}) \times R^{3/2}
\subset A(R, R^2). 
$$
Now we show the following claim:

\parag{\bf Claim.}
{\it
There exists a quasiconformal homeomorphism
$$
\Theta_n^0 : \overline{A(R, R^2)} \to 
              \overline{U}_{s_n} \smallsetminus V_{s_n}'
$$ 
for $s_n$ such that 
\begin{itemize}
\setlength{\itemsep}{0.5mm} 
\item
$\Theta_n^0$ is quasiconformal;

\item
$\Theta_n^0$ is equivariant on the boundary, i.e.,
$\Theta_n^0(Z^2) = G_{s_n}(\Theta_n^0(Z))$ for $|Z| = R$;

\item
$\Theta_n^0(Z) = \Psi_{s_n}^{-1}(R^{- 3/2}Z)$ 
for $Z \in \Gamma_0(c_0+\eta_n)$; and

\item
$\Theta_n^0(\Gamma_0(c_0+\eta_n)) = J(f_{s_n})$. 
\end{itemize}
}

\parag{Proof of the Claim.}
Let $U$, $U'$, and $V'$ be the quasidisks that are the images of 
$U_{s_n}$, $U_{s_n}'$, and $V_{s_n}'$ by the straightening map 
$\Psi_{s_n}$ respectively. 
Then $P_{c_0 +\eta_n}:U' \to U$ is a quadratic-like restriction 
of $P_{c_0 +\eta_n}$.
We also take a smooth Jordan domain $V''$ such that
$V' \Subset V'' \Subset U'$
and 
$J(P_{c_0+\eta_n}) \subset U' \sminus \overline{V''}$.

Let $R>1$ be large enough to 
have $D(R^{-1/2}) \Subset V''$ and $U' \Subset D(R^{1/2})$.
We first show that there exists a quasiconformal map
$\psi:\Bar{A(R^{-1/2}, R^{1/2})}
 \to \Bar{U} \sminus V'$
that fixes the closed annulus $\Bar{U'} \sminus V''$
including the Julia set $J(P_{c_0+\eta_n})$.
To do this, we construct three quasiconformal maps
$\psi_0:\Bar{D(R^{1/2})} \sminus U' \to \Bar{U} \sminus U'$,
$\psi_1=\mathrm{id}: 
\Bar{U'} \sminus V'' \to \Bar{U'} \sminus V''$,
and 
$\psi_2:\Bar{V''} \sminus D(R^{-1/2})
\to \Bar{V''} \sminus V'$
and glue them together along the boundaries.
Indeed, we can find such a $\psi_0$ because 
there are conformal isomorphisms
from 
$D(R^{1/2}) \sminus \overline{U'}$
and 
$U \sminus \overline{U'}$ 
to round annuli of finite modulus,
and composition with a radial stretch from one of the round annuli
to the other will give a quasiconformal map 
from $D(R^{1/2}) \sminus \overline{U'}$ to $U \sminus \overline{U'}$.
Since the boundary component of these annuli are 
quasicircles, $\psi_0$ extends homeomorphically 
to the boundary and we can modify the map  
in the category of quasiconformal maps  
in such a way that $\psi_0$ agrees with $\psi_1=\mathrm{id}$ 
on the inner boundary $\partial U'$. 
By the same reasoning we can find $\psi_2$ that 
agrees with $\psi_1=\mathrm{id}$ on $\partial V''$.

Next we define $\widetilde{\psi}:\overline{A(R,R^2)} \to \Bar{U_{s_n}} \sminus V_{s_n}$ by
$\widetilde{\psi}(Z):=\Psi_{s_n}^{-1} \cc \psi(R^{-3/2} Z)$ 
for each $Z \in \overline{A(R,R^2)}$.
This is a quasiconformal map sending $\Gamma_0(c_0+\eta_n)$
to the Julia set $J(f_{s_n})$.

Finally, to obtain $\Theta_n^0$ in the statement, 
we modify $\widetilde{\psi}$ near the inner boundary $|Z|=R$
to satisfy $\Theta_n^0(Z^2) 
= G_{s_n}(\Theta_n^0(Z))$ for $|Z| = R$.
This is again possible because 
$\partial V_{s_n}'$ and $\partial V'$ are quasicircles. 
\hfill\QED{\small (Claim)}\\

Let us return to the proof of Lemma \ref{lem_G_is_M_like}.
The Julia set
$J(f_{c}) \subset U_c' \smallsetminus \overline{V_c'}$ is a Cantor set
for every $c \in W_n$ for the same reason for $J(f_{s_n})$ and this, 
as well as $\partial U_c$ and $\partial V_c'$ undergo holomorphic motion 
(see \cite[p.229]{Shishikura 1998}).
By S{\l}odkowski's theorem (\cite{Slodkowski 1991}) there exists a 
holomorphic motion $\iota_c$ on $\C$
which induces these motions. Finally 
define $\Theta_c := \iota_c \circ \Theta_n^0$ , then 
$\Theta = \Theta_n := \{ \Theta_c \}_{c \,\in\, W_n}$ is a desired
tubing for $\boldsymbol{G}_n$. To prove this, 
we have to check that
$\boldsymbol{G}_n$ with $\Theta_n$
satisfies the conditions (1)--(8) for a Mandelbrot-like family.
The condition (1) is already shown in Lemma~\ref{lem_4.2}. 
It is easy to check the conditions (2)--(7). 
Finally the one turn condition (8) is proved as follows: 
Note that $\check{c}_n(\gamma(t))$ satisfies
$$
  \tau(\check{c}_n(\gamma(t))) + n = v(\check{c}_n(\gamma(t)), \gamma(t))
$$
When $c$ ranges over $\partial W_n$ making one turn, the variable $t$ for
both sides varies from $t=0$ to $t=1$. 
Since $v(\check{c}_n(\gamma(t)),\gamma(t))$ 
is very close to $v(c_1, \gamma(t))$, which is a parameterization of
$\partial \wt{V}_{c_1}$ for sufficiently large $n$, 
$v(\check{c}_n(\gamma(t)),\gamma(t))$ and hence
$\tau(\check{c}_n(\gamma(t))) + n$ 
makes one turn in a very thin tubular neighborhood of 
$\partial \wt{V}_{c_1}$ as $t$ moves from 0 to 1. 
This implies that 
$f_c^{k\nu n}(0) = f_c^L(0)$ 
makes one turn in a very thin tubular neighborhood of $\partial V_{c_1}$. 
Hence $G_c(0)-0 = G_c(0) = g_c \circ f_c^{k\nu n}(0) = g_c \circ f_c^L(0)$ makes
one turn in a very thin tubular neighborhood of $\partial U_{c_1}$. 
In particular this shows that $G_c(0)-0$ makes one turn around $0 \in U_{c_1}$.
\hfill \QED(Lemma \ref{lem_G_is_M_like})

\section{
Step (P4): Existence of a copy of $\cM(c_0+\eta)$ in $W$
}

Finally we show that the family 
$\bs{G}$ provides the desired quasiconformal
copy of $\cM(c_0+\eta)$ in $W$.

\parag{End of the proof of Theorem \ref{thm_A'}.}
Recall that in Step (P1),
we took the sector $S$ such that 
$S \subset D(c_1, \e')$ for the given $\e'>0$ in the statement.
Hence $W=W_n \subset S$ is contained in $D(c_1, \e')$.
Its center $s=s_n$ tends to $c_1$ as $n \to \infty$
by construction. Hence for any $\vep > 0$
given in the statement, 
we may assume that 
$|\chi_{s_0}(s)-\chi_{s_0}(c_1)|<\e$
by taking a sufficiently large $n$.
(Here we used the continuity of 
$\chi_{s_0}:\Lambda \to \chi_{s_0}(\Lambda)$.)
Set $\eta:=\chi_{s_0}(s)-\chi_{s_0}(c_1)=\sigma-c_0$
such that $\sigma=c_0+\eta \in \C \sminus M$.
Now we will show that $\cM(\sigma)$ appears 
quasiconformally in $M \cap \Bar{W}$.

Let $\Theta=\{\Theta_c\}_{c\, \in \, W}$ 
be the tubing of $\bs{G}$ given in Lemma \ref{lem_G_is_M_like}.
Let us check that the set
$$
{\mathcal N}
:=
M_{\boldsymbol G}
\cup 
\{ c \ | \ 
     G_c^k(0) \in \Theta_c(\Gamma_0(\sigma)) \quad \text{for some} \ k \in \N \}
$$
is the image of ${\mathcal M}(\sigma)$ 
by the quasiconformal map 
$\chi_{\Theta}^{-1}:W_R \to W$, 
where $M_{\boldsymbol G}=\chi_{\Theta}^{-1}(M)$ 
is the connectedness locus of ${\boldsymbol G}$. 
Indeed, by \eqref{eq_chi_Theta},
$G_c^k(0) \in \Theta_c(\Gamma_0(\sigma))$ 
is equivalent to 
$\braces{\Phi_M(\chi_\Theta(c))}^{2^{k-1}} \in \Gamma_0(\sigma)$
and thus $\Phi_M(\chi_\Theta(c)) \in  \Gamma_{k-1}(\sigma)$.
This implies that $\chi_\Theta(c) \in \cM(\sigma)$.

Finally we check that ${\mathcal N} \subset M$ and 
$\partial {\mathcal N} \subset \partial M$.
If $c \in M_{\bs{G}}$, 
then the orbit of the critical point $0$ by 
$G_c = g_c \circ f_c^{k\nu n} = P_c^{k\nu np + N}$ is bounded.
Hence we have $c \in M$.
If $c$ satisfies 
$$
G_c^k(0) \in \Theta_c(\Gamma_0(\sigma))=J(f_c)
$$ 
for some $k \in \N$,
it implies $c \in M$ as well 
since $J(f_c)$ is invariant under $f_c=P_c^p$.
So the set $\mathcal N$ is a subset of $M$.
To show $\partial {\mathcal N} \subset \partial M$,
it is enough to show that 
the Misiurewicz parameters 
are dense in $\partial {\mathcal N}$.
They are clearly dense in $\partial M_{\bs{G}}$.
Let $\mathcal{L}$ be the set of $c$ such that $G_c^k(0)$ is a 
repelling periodic point in $J(f_c)$ for some $k \ge 1$.
Then $\mathcal{L}$ is a dense subset of
$\partial {\mathcal N} \sminus \partial M_{\bs{G}} 
  ={\mathcal N} \sminus  M_{\bs{G}}$.
This completes the proof of Theorem \ref{thm_A'}.
\QED

\parag{Acknowledgment.} 
The authors would like to thank the referee for careful reading and useful comments. 
They were partly supported by JSPS KAKENHI Grants 
16K05193, 17K05296, and 19K03535.

\end{document}